\documentclass[a4paper,11pt]{article}
\usepackage[latin1]{inputenc}
\usepackage[francais,english]{babel}    
\usepackage{amssymb}
\usepackage{amsmath}
\usepackage{amsthm}
\usepackage{multicol}
\usepackage{graphicx}
\usepackage{color}
\usepackage[T1]{fontenc}
\usepackage{yfonts}
\usepackage{placeins}
\usepackage{array}
\usepackage{dsfont}
\usepackage{stmaryrd}

\newcommand{\sign}{\begin{flushright}
Thomas Haettel \\
Université Paris-Sud 11 \\
Département de Mathématiques \\
UMR 8628 CNRS \\
91405 Orsay \\
thomas.haettel@math.u-psud.fr
\end{flushright}}

\newtheorem{thm}{Théorème}[section]
\newcommand{\bthm}{\begin{thm}}
\newcommand{\ethm}{\end{thm}}

\newtheorem*{thn}{Théorème}
\newcommand{\bthn}{\begin{thn}}
\newcommand{\ethn}{\end{thn}}

\newtheorem{defi}[thm]{Définition}
\newcommand{\bdf}{\begin{defi}}
\newcommand{\edf}{\end{defi}}

\newtheorem{qc}{Question de Cours}
\newcommand{\bqc}{\begin{qc}}
\newcommand{\eqc}{\end{qc}}

\newtheorem{rqc}{R\'{e}ponse}
\newcommand{\brqc}{\begin{rqc}}
\newcommand{\erqc}{\end{rqc}}

\newtheorem{exo}{Exercice}
\newcommand{\bex}{\begin{exo}}
\newcommand{\eex}{\end{exo}}

\newtheorem{sol}{Solution}
\newcommand{\bsol}{\begin{sol}}
\newcommand{\esol}{\end{sol}}

\newtheorem{pro}[thm]{Proposition}
\newcommand{\bpro}{\begin{pro}}
\newcommand{\epro}{\end{pro}}

\newtheorem{cor}[thm]{Corollaire}
\newcommand{\bcor}{\begin{cor}}
\newcommand{\ecor}{\end{cor}}

\newtheorem{lem}[thm]{Lemme}
\newcommand{\blem}{\begin{lem}}
\newcommand{\elem}{\end{lem}}

\newtheorem*{rmq}{Remarque}
\newcommand{\brq}{\begin{rmq} \upshape}
\newcommand{\erq}{\end{rmq}}

\newtheorem*{exe}{Exemple}
\newcommand{\bexe}{\begin{exe} \upshape}
\newcommand{\eexe}{\end{exe}}

\newtheorem*{pre}{Démonstration}
\newcommand{\bp}{\begin{pre} \upshape}
\newcommand{\ep}{\hfill \qed \end{pre}}
\newcommand{\epp}{\end{pre}}

\newcommand{\beq}{\begin{eqnarray*}}
\newcommand{\eeq}{\end{eqnarray*}}

\newcommand{\beqn}{\begin{equation}}
\newcommand{\eeqn}{\end{equation}}

\newcommand{\ben}{\begin{enumerate}}
\newcommand{\een}{\end{enumerate}}
\newcommand{\bit}{\begin{itemize} \renewcommand{\labelitemi}{$\bullet$} }
\newcommand{\eit}{\end{itemize}}

\newcommand{\bfg}{
\begin{figure}[H]
\begin{center}}
\newcommand{\efg}{
\end{center}
\end{figure}
\FloatBarrier}

\newcolumntype{M}[1]{>{\raggedright}m{#1}}

\newcommand{\df}{\emph}

\newcommand{\R}{\mathbb{R}}
\newcommand{\Q}{\mathbb{Q}}
\newcommand{\N}{\mathbb{N}}
\newcommand{\Z}{\mathbb{Z}}
\newcommand{\C}{\mathbb{C}}

\renewcommand{\P}{\mathbb{P}}
\newcommand{\K}{\mathbb{K}}

\renewcommand{\H}{\mathbb{H}}
\newcommand{\U}{\operatorname{U}}

\newcommand{\bs}{\symbol{92}}

\newcommand{\ov}{\overline}
\newcommand{\OO}{\mathcal{O}}
\newcommand{\M}{\mathcal{M}}

\newcommand{\sh}{\operatorname{sh}}

\newcommand{\argsh}{\operatorname{argsh}}

\renewcommand{\Im}{\operatorname{Im}}
\renewcommand{\Re}{\operatorname{Re}}
\newcommand{\Ker}{\operatorname{Ker}}

\renewcommand{\t}{ ^t\!}
\newcommand{\tr}{\operatorname{tr}}

\newcommand{\Diag}{\operatorname{Diag}}

\newcommand{\Card}{\operatorname{Card}}

\newcommand{\Id}{\operatorname{Id}}

\newcommand{\Sp}{\operatorname{Sp}}

\newcommand{\End}{\operatorname{End}}
\newcommand{\Tor}{\operatorname{Tor}}
\newcommand{\Teich}{\operatorname{Teich}}
\newcommand{\Mod}{\operatorname{Mod}}
\newcommand{\MCG}{\operatorname{MCG}}
\newcommand{\sns}{\operatorname{sns}}

\newcommand{\lb}{\llbracket}
\newcommand{\rb}{\rrbracket}
\newcommand{\ra}{\rightarrow}
\newcommand{\ral}[1]{\underset{#1}{\longrightarrow}}

\renewcommand{\geq}{\geqslant}
\renewcommand{\leq}{\leqslant}

\newcommand{\gothique}{\mathfrak}

\newcommand{\hh}{\gothique{h}}

\renewcommand{\log}{\operatorname{log}}

\newcommand{\id}{\operatorname{id}}

\newcommand{\Aut}{\operatorname{Aut }}
\newcommand{\Int}{\operatorname{Int}}
\newcommand{\GL}{\operatorname{GL}}
\newcommand{\SL}{\operatorname{SL}}
\newcommand{\SO}{\operatorname{SO}}
\renewcommand{\O}{\operatorname{O}}
\newcommand{\SU}{\operatorname{SU}}

\newcommand{\Spin}{\operatorname{Spin}}
\newcommand{\Sym}{\operatorname{Sym}}
\newcommand{\Diff}{\operatorname{Diff}}

\newcommand{\<}{\langle}
\renewcommand{\>}{\rangle}

\makeatletter
\def\Ddots{\mathinner{\mkern1mu\raise\p@
\vbox{\kern7\p@\hbox{.}}\mkern2mu
\raise4\p@\hbox{.}\mkern2mu\raise7\p@\hbox{.}\mkern1mu}}
\makeatother

\makeatletter
\def\maketitles{%
  \null
  \thispagestyle{empty}%
  \vfill
  \begin{center}\leavevmode
    \normalfont
    {\LARGE \@title\par}%
    \vskip 1.2cm
    {\large \@author\par}%
    \vskip 1.2cm
    {\large \@subtitle\par}%
    \vskip 0.8cm
    {\large \@date\par}%
  \end{center}%
  \vfill
  \null
  \cleardoublepage
  }

\def\date#1{\def\@date{#1}}
\def\author#1{\def\@author{#1}}
\def\title#1{\def\@title{#1}}
\def\subtitle#1{\def\@subtitle{#1}}
\makeatother

\title{Compactification de Thurston d'espaces de réseaux marqués et de l'espace de Torelli}
\author{Thomas Haettel}
\date{\today}

\begin{document}

\selectlanguage{francais}

\maketitle

\selectlanguage{francais}

\selectlanguage{english}
\begin{abstract} We define a compactification of symmetric spaces of noncompact type, seen as spaces of isometry classes of marked lattices, analogous to the Thurston compactification of the Teichmüller space, and we show that it is equivariantly isomorphic to a Satake compactification. We then use it to define a new compactification of the Torelli space of a hyperbolic surface with marked points, and we show that it is equivariantly isomorphic to the Satake compactification of the image of the period mapping. Finally, we describe the natural stratification of a subset of the boundary.
\end{abstract}
\selectlanguage{francais}
\section*{Introduction}

Soit $S$ une surface lisse compacte connexe orientable, de genre $g \geq 2$. La compactification de Thurston de l'espace de Teichmüller $\Teich(S)$ de $S$ est l'adhérence du plongement $\Teich(S) \ra \P({\R_+}^{\!\Pi_1(S)})$ défini par les distances de translation. Cette compactification est homéomorphe à la boule fermée de dimension $6g-6$, munie d'une action continue du groupe modulaire $\MCG(S)$, et elle a permis à Thurston de classer à isotopie près les difféomorphismes de $S$ (voir par exemple~\cite{thurston}).

Tout espace symétrique de type non compact classique peut être vu comme un espace de réseaux marqués d'un espace euclidien ou hermitien (voir par exemple~\cite{bavard}). Soit $m \geq 1$ un entier, on appelle réseau marqué de covolume $1$ de $\R^m$ tout morphisme injectif de $\Z$-modules $f:\Z^m \ra \R^m$ dont l'image est discrète et de covolume $1$. L'espace des classes d'isométrie de tels réseaux marqués est naturellement homéomorphe à l'espace symétrique de type non compact $\mathcal{E}_m = \SL_m(\R) / \SO_m(\R)$. On peut définir une compactification de Thurston de cet espace, grâce aux distances de translation (où $\R^m$ est muni de la norme euclidienne usuelle) :
\beq \phi : \mathcal{E}_m & \ra & \P({\R_+}^{\!\Z^m}) \\
\ [f:\Z^m \ra \R^m] & \mapsto & [u \mapsto ||f(u)||\,].\eeq

\bthn L'application $\phi$ est un plongement : la \df{compactification de Thurston} de l'espace symétrique $\mathcal{E}_m$ est l'adhérence $\ov{\mathcal{E}_m}^\mathcal{T}$ de l'image de $\phi$. Cette compactification est isomorphe (de manière $\SL_m(\Z)$-équivariante) à la compactification de Satake $\ov{\mathcal{E}_m}^\mathcal{S}$ associée à la représentation tautologique de $\SL_m(\R)$ sur $\R^m$. \ethn

Ceci répond à une question posée par Frédéric Paulin lors d'un exposé Bourbaki (voir~\cite{paulin_bourbaki}). Nous démontrons en fait un résultat plus général (théorème~\ref{thm:thurston_satake}), pour le corps $\R$, $\C$ ou le corps gauche $\H$ des quaternions de Hamilton, ainsi que pour des sous-espaces de réseaux autoduaux, ce qui permet de traiter le cas de tous les espaces symétriques de type non compact classiques. Et nous étendons ce théorème au cas du groupe de Lie exceptionnel $E_{6(-26)}$ (théorème~\ref{thm:thurston_satake_O}), qui est la forme réelle non compacte de rang réel $2$ du groupe de Lie complexe exceptionnel $E_6$.

\bigskip

Puis nous appliquons cette construction pour définir une compactification à la Thurston de l'espace de Torelli. L'espace de Torelli $\Tor(S)$ de la surface $S$ est l'espace des classes d'isotopie de surfaces hyperboliques munies d'un marquage de leur cohomologie par celle de $S$. L'espace de Torelli $\Tor(S)$ est le quotient de l'espace de Teichmüller par le groupe de Torelli $T(S)$, qui est le sous-groupe du groupe modulaire de $S$ constitué des classes d'isotopie de difféomorphismes de $S$ ayant une action triviale en homologie. Un théorème de Mess (voir~\cite{mess}) énonce qu'en genre $2$, le groupe de Torelli est un groupe libre sur une infinité dénombrable de générateurs. En genre supérieur ou égal à $3$, le groupe de Torelli est de type fini (voir~\cite{johnson}), mais l'une des grandes questions est de savoir s'il est de présentation finie (voir~\cite{mcg}).

Dans l'espoir de mieux comprendre le type d'homotopie de l'espace de Torelli, nous allons en définir une compactification naturelle. Considérons une surface hyperbolique $X$ marquée par $h : S \ra X$. Le théorème de Hodge identifie l'espace $H^1(X,\R)$ avec l'espace des $1$-formes différentielles harmoniques sur $X$, espace qui est muni du produit scalaire $L^2$ : notons $\|\cdot\|_X$ la norme euclidienne associée sur $H^1(X,\R)$.
Considérons alors l'application
\beq \psi : \Tor(S) & \ra & {\R_+}^{\!H^1(S,\Z)} \\
\ [X,h] & \mapsto & \left\{ \omega \mapsto  \| h^*(\omega) \|_X \right\}, \eeq
et $\ov{\psi} : \Tor(S) \ra \P({\R_+}^{\!H^1(S,\Z)})$ son application quotient. Nous montrons dans la partie~\ref{sec:torelli} que l'application $\ov{\psi}$ est un revêtement de degré $2$ sur son image, ramifié sur le lieu hyperelliptique de $\Tor(S)$. En considérant l'adhérence de l'image de $\ov{\psi}$ dans l'espace compact $\P({\R_+}^{\!H^1(S,\Z)})$, nous définissons une compactification de l'espace de Torelli, que nous appelons \df{compactification de Thurston} de $\Tor(S)$.

Nous allons comparer notre compactification à une autre compactification naturelle de l'espace de Torelli. L'application période envoie l'espace de Torelli $\Tor(S)$ dans l'espace symétrique $\Sp_{2g}(\R) / \SU(g)$. L'adhérence de l'image de l'application période dans la compactification de Satake de $\Sp_{2g}(\R) / \SU(g)$ associée à la représentation standard de $\Sp_{2g}(\R)$ sur $\R^{2g}$ définit une compactification de l'espace de Torelli, que nous appelerons \df{compactification de Satake} de $\Tor(S)$.

\bthn  Les compactifications de Thurston et de Satake de l'espace de Torelli sont isomorphes de manière $\Sp_{2g}(\Z)$-équivariante. \ethn

De plus, nous décrivons une partie de bord de cette compactification : nous allons décrire l'adhérence de l'image de $\psi$ (et il serait intéressant d'avoir une description analogue, au moins à homotopie près, pour $\ov{\psi}$).

Considérons $K^{sep}(S)$ le complexe des courbes séparantes de $S$ : les sommets de ce complexe simplicial sont les classes d'homotopie de courbes fermées simples séparantes non triviales, et les ($k$-$1$)-simplexes sont les $k$-uplets $\sigma = \{[\gamma_1], \ldots , [\gamma_k]\}$ de classes de telles courbes deux à deux disjointes et non homotopes. Notons également $\Sigma K^{sep}(S)$ l'ensemble des simplexes de $K^{sep}(S)$.

Si $\sigma$ est un tel ($k$-$1$)-simplexe, considérons les $k$+$1$ composantes connexes de $S \bs \cup_{j=1}^k \gamma_j$ , et fixons pour chacune d'elles un homéomorphisme avec $S_j \bs P_j$, où $S_j$ est une surface compacte lisse sans bord de genre $g_j \geq 1$, et où $P_j \subset S_j$ est un ensemble fini de points. Notons
\beq \psi_\sigma : \Tor_\sigma(S) = \prod_{j=0}^k \Tor(S_j) & \ra & {\R_+}^{\!H^1(S,\Z)} \\
([X_j,h_j])_{0 \leq j \leq k} & \mapsto & \left\{ \omega = \sum_{j=0}^k \kappa_j^* \omega_j \mapsto \left( \sum_{j=0}^k \| {h_j^{-1}}^*(\omega_j) \|^2 \right)^{\frac{1}{2}} \right\},\eeq
où $\omega_j \in H^1(S_j,\Z)$, et où $\kappa_j : S \ra S_j$ est l'application d'écrasement des $S_{j'}$, pour $j' \neq j$ (voir partie~\ref{subsec:torelli_strates}).

\bthn L'adhérence de l'image de l'application $\psi$ dans l'espace ${\R_+}^{\!H^1(S,\Z)}$ est la réunion disjointe des strates
$$ \ov{\psi(\Tor(S))} = \psi(\Tor(S)) \sqcup \bigsqcup_{\sigma \in \Sigma K^{sep}(S)} \psi_\sigma(\Tor_\sigma(S)).$$ \ethn

Dans une première partie, nous rappelons quelques résultats élémentaires d'algèbre linéaire quaternionique (voir~\cite{bourbaki_algebre9}). Puis nous définissons la compactification de Thurston de l'espace des réseaux marqués d'un espace euclidien, hermitien complexe ou quaternionique. Ensuite, nous établissons l'isomorphisme avec la compactification de Satake. Dans une quatrième partie, nous étendons ce résultat aux espaces de réseaux autoduaux.

Ensuite, nous montrons que ces résultats s'appliquent à l'espace de Torelli. Et enfin, nous décrivons la stratification d'une partie du bord de cette compactification.

\bigskip

Je tiens à remercier chaleureusement Frédéric Paulin pour les nombreux conseils qu'il m'a donnés sur cet article. Je tiens également à remercier Daniel Massart, qui m'a expliqué les liens entre les normes stable et $L^2$.

\section{Un peu d'algèbre linéaire quaternionique}

\label{sec:algebre_lineaire}

Soit $\K$ le corps (commutatif) $\R$, $\C$ ou le corps (gauche) des quaternions de Hamilton $\H$ (de base vectorielle réelle $1,i,j,k$, où $i^2=j^2=-1$ et $ij=-ji=k$). Si $x \in \K$, notons $\ov{x}$ son conjugué et $\tr(x)$ sa trace réduite (avec $\ov{x}=x$ si $\K=\R$, $\ov{x}=a-ib$ si $x=a+ib \in \K=\C$ et $\ov{x}=a-ib-jc-kd$ si $x=a+ib+jc+kd \in \K=\H$). Nous avons des inclusions évidentes $\R \subset \C \subset \H$. Pour la théorie générale des formes sesquilinéaires sur des espaces vectoriels de dimension finie sur des corps gauches, nous renvoyons à~\cite{bourbaki_algebre9}.

Fixons un entier $m \geq 1$. L'espace vectoriel $\K^m$ sera toujours considéré à droite, et l'action linéaire des matrices de $\M_m(\K)$ sur $\K^m$ à gauche. Munissons l'espace vectoriel $\K^m$ de sa structure hermitienne standard $\<x,y\> = \sum_{i=1}^m \ov{x_i}y_i$.

Si $M \in \M_m(\K)$, notons $\ov{M}$ la matrice de coefficients les conjugués de ceux de $M$, et $M^* = \,\t\ov{M}$ l'adjoint de $M$ : il vérifie $(MN)^* = N^*M^*$. Une matrice $M$ est dite \df{hermitienne} si $M^*=M$, et \df{unitaire} si $M M^*=I_m$ (ce qui équivaut à $M^* M=I_m$), où $I_m$ désigne la matrice identité de $\M_m(\K)$.

Rappelons une définition du déterminant de Dieudonné d'une matrice à coefficients dans $\H$ (voir par exemple~\cite[§~IV.1, p.~149]{artin}). Pour cela, considérons l'application $\eta$ de $\M_m(\H)$ dans $\M_{2m}(\C)$ qui à une matrice $M=A+jB$, où $A,B \in M_m(\C)$ associe la matrice par blocs
$$ \left( \begin{array}{cc} A & -\ov{B} \\ B & \ov{A} \end{array} \right) \in \M_{2m}(\C).$$
C'est un homomorphisme injectif de $\R$-algèbres à droite, équivariant pour l'adjoint. Remarquons de plus que le déterminant de $\eta(M)$ est un nombre réel positif. Ceci permet de définir $\det(M) = \sqrt{\det(\eta(M))}$, le \df{déterminant de Dieudonné} de $M$. La matrice $M$ est alors inversible si et seulement si $\det(M) \neq 0$, et si $N \in \M_m(\H)$ alors $\det(MN)=\det(M)\det(N)$.
Cela permet de définir le sous-groupe $\SL_m(\H)$ de $\GL_m(\H)$. Ce sous-groupe peut également être défini intrinsèquement : c'est le sous-groupe de $\GL_m(\H)$ constitué des automorphismes linéaires de $\H^m$ qui préservent une mesure de Haar de $\H^m$.

Notons $\U_m(\K)$ le groupe unitaire de $\K^m$ : c'est le sous-groupe de $\GL_m(\K)$ constitué des matrices unitaires. Notons de plus $\SU_m(\K) = \U_m(\K)~\cap~\SL_m(\K)$ le groupe spécial unitaire de $\K^m$. Remarquons que lorsque $\K=\H$, le groupe $\SU_m(\H)$ est égal à $\U_m(\H)$, et également à l'image réciproque par l'application $\eta$ du sous-groupe $U_{2m}(\C)$. Le groupe $\SU_m(\K)$ est un sous-groupe compact maximal du groupe de Lie $\SL_m(\K)$.

\blem[Diagonalisation des matrices hermitiennes] \label{lem:diag} Soit $M \in \M_m(\H)$ une matrice hermitienne positive. Alors il existe une matrice unitaire $U \in \U_m(\H)$ et une matrice diagonale positive réelle $D \in \M_m(\R)$ telles que $M = UDU^{-1}$. \elem

\bp Pour le théorème général de réduction des matrices hermitiennes, voir~\cite[Théorème~1, p.~90]{bourbaki_algebre9}. Nous donnons ici une preuve rapide en admettant la diagonalisation des matrices hermitiennes complexes. L'application $\eta : \M_m(\H) \ra \M_{2m}(\C)$ a pour image l'ensemble des matrices des endomorphismes qui commutent avec l'endomorphisme réel $\alpha$ de $\C^{2m}$ défini par $X \mapsto J\,\ov{X}$, où $J$ est la matrice
$$ J = \left( \begin{array}{cc} 0 & I_m \\ -I_m & 0 \end{array} \right).$$
Soit $N=\eta(M) \in \M_{2m}(\C)$ : puisque l'application $\eta$ est équivariante pour l'adjoint, la matrice $N$ est hermitienne positive. Les sous-espaces propres complexes de $N$ sont orthogonaux et stables par l'endomorphisme réel $\alpha$, donc il existe une matrice unitaire $U' \in \U_{2m}(\C)$ commutant avec $\alpha$ telle que $D'=U'^{-1}NU'$ soit diagonale positive réelle. Puisque $U'$ et $D'$ commutent avec $\alpha$, ils appartiennent à l'image de $\eta$ : soit donc $U \in \SU_m(\H)$ et $D$ diagonale positive réelle tels que $U'=\eta(U)$ et $D'=\eta(D)$. Ainsi $M = UDU^{-1}$. \ep

\blem[Décomposition polaire] \label{lem:polaire} Soit $M \in \M_m(\K)$. Alors il existe une matrice unitaire $U \in \U_m(\K)$ et une unique matrice hermitienne positive $P \in \M_m(\K)$ telles que $M = PU$. \elem

\bp Montrons tout d'abord ce résultat si $M$ est inversible. La matrice $MM^*$ est hermitienne positive, donc d'après la diagonalisation des matrices hermitiennes positives (lemme~\ref{lem:diag} si $\K=\H$), considérons $P$ l'unique matrice hermitienne positive telle que $P^2 = MM^*$. La matrice $P$ est inversible, et la matrice $U=P^{-1}M$ est alors unitaire.

Si $M$ n'est pas inversible, par densité de $\GL_m(\K)$ dans $\M_m(\K)$, considérons une suite $(M_n)_{n \in \N}$ de matrices inversibles convergeant vers $M$. Pour tout $n \in \N$, considérons une décomposition polaire $M_n = P_nU_n$ de la matrice $M_n$. Par compacité du sous-groupe $U_m(\K)$, on peut supposer quitte à extraire que la suite $(U_n)_{n \in \N}$ converge vers une matrice unitaire $U$. Dans ce cas, la suite $(P_n = M_n {U_n}^{-1})_{n \in \N}$ converge vers une matrice $P=MU^{-1}$ hermitienne positive, et nous avons $M=PU$. Et la matrice $P$ est l'unique matrice hermitienne positive telle que $P^2 = MM^*$. 
\ep

Nous aurons également besoin d'identités de polarisation, la première étant immédiate.

\blem[Identité de polarisation complexe] \label{lem:polarisation_C} Soit $\tau=a+ib \in \C \bs \R$. Alors pour tous $u,v \in \C^m$, nous avons :
$$ \<u|v\> = \frac{1}{4} \left( (1+\frac{ia}{b})\left( \| u+v \| ^2 -  \| u-v \| ^2\right) - \frac{i}{b}\left( \| u+v \tau  \| ^2 -  \| u-v \tau  \| ^2\right) \right).$$
\hfill $\qed$ 
\elem

\blem[Identité de polarisation quaternionique] \label{lem:polarisation_H} Soit $(q_1,\ldots,q_4)$ une $\R$-base de $\H$. Alors il existe quatre quaternions $(\lambda_1,\ldots,\lambda_4) \in \H^8$ tels que, pour tous $u,v \in \H^m$, nous ayons :
$$ \<u|v\> = \sum_{l=1}^4 \left(  \| u+ vq_l \| ^2 -  \| u- vq_l \| ^2\right) \lambda_l .$$
\elem

\bp Remarquons tout d'abord que
$$\left(  \| u+ v q_l \| ^2 -  \| u- v q_l \| ^2\right) = 4 \Re(\<u, v q_l\>) = 4 \Re(\ov{u}vq_l) = 4 \Re(\<u,v\>q_l).$$
Considérons l'application $\R$-linéaire
\beq \varphi : \H & \ra & \R^4 \\
q & \mapsto & (4\Re(qq_l))_{l \in \lb 1,4 \rb}.\eeq
Si $q \in \H$ appartient au noyau de $\varphi$, alors par $\R$-linéarité $4\Re(q\ov{q}) = 4|q|^2 = 0$ donc $q=0$ : l'application $\varphi$ est donc un isomorphisme $\R$-linéaire. Pour tout $l \in \lb 1,4 \rb$, notons $\lambda_l = \varphi^{-1}(f_l)$, où $(f_1,\ldots,f_4)$ désigne la base canonique de $\R^4$.
Alors, pour tous $u,v \in \H^m$, nous avons par $\R$-linéarité
\beq \<u,v\> &=& \varphi^{-1}(\varphi(\<u,v\>)) = \varphi^{-1}\left((4\Re(\<u,v\> q_l))_{l \in \lb 1,4 \rb}\right) \\
&=& \varphi^{-1}\left(\sum_{l=1}^4 (4\Re(\<u,v\> q_l)) f_l \right) \\
&=& \sum_{l=1}^4 4\Re(\<u,v\>) \varphi^{-1}(f_l) \\
&=& \sum_{l=1}^4 4\Re(\<u,v\>) \left(  \| u+ v q_l \| ^2 -  \| u- v q_l \| ^2\right) \lambda_l. \eeq
\ep

Nous aurons également besoin du lemme suivant.

\blem \label{lem:norme} Soient $M$ et $M'$ deux matrices de $\M_m(\K)$ telles que, pour tout $u \in \K^m$, nous ayons $ \| M(u) \| = \| M'(u) \| $. Alors il existe $K \in \U_m(\K)$ tel que $KM=M'$. \elem

\bp D'après la décomposition polaire (lemme~\ref{lem:polaire} si $\K=\H$) et modulo multiplication à gauche par deux éléments de $U_m(\K)$, on peut supposer que $M$ et $M'$ sont hermitiennes positives. D'après la réduction des matrices hermitiennes (lemme~\ref{lem:diag}) et modulo multiplication à gauche par deux éléments de $U_m(\K)$, on peut supposer qu'il existe deux matrices diagonales positives $D,D'$ et deux matrices unitaires $U,U'$ telles que $M=DU$ et $M'=D'U'$. Alors, pour tout $u \in \K^m$, nous avons $ \| Du \|  =  \| D'U'U^{-1}u \| $. Ainsi, si $u$ appartient au sous-espace propre de $D$ associé à la valeur propre maximale, on en déduit que $U'U^{-1}u$ appartient au sous-espace propre de $D'$ associé à la valeur propre maximale. Puis on montre que si $u$ appartient au sous-espace propre de $D$ associé à la deuxième plus grande valeur propre (s'il y en a), alors $U'U^{-1}u$ appartient au sous-espace propre de $D'$ associé à la deuxième plus grande valeur propre. Finalement, on montre que la matrice unitaire $U'U^{-1}$ envoie les sous-espaces propres de $D$ sur ceux de $D'$. Quitte à permuter les coefficients de $D$ ou $D'$, on en déduit que $D=D'$ et $D(U'U^{-1}) = (U'U^{-1})D$. D'où $M'= DU' = D(U'U^{-1})U = (U'U^{-1})DU = KM$, avec $K = U'U^{-1} \in \U_m(\K)$. \ep

\section{Compactification de l'espace des réseaux hermitiens marqués}

\label{sec:thurston}

Si $\K=\R$, notons $\OO=\Z$. Si $\K=\C$, notons $\OO$ un ordre de l'anneau des entiers d'un corps de nombres quadratique imaginaire (par exemple $\OO=\Z[i]$). Enfin si $\K=\H$, notons $\OO$ un ordre dans une algèbre de quaternions $A$ sur $\Q$ non déployée sur $\R$ --- et nous identifierons alors $A \otimes_\Q \R$ et $\H$ --- (par exemple $\OO=\Z[1,i,j,\frac{1+i+j+k}{2}]$). Ainsi $\OO$ est un $\Z$-réseau de l'espace vectoriel réel $\K$, et est un anneau contenant $1$.

On appelle \df{réseau} (ou $\OO$-réseau lorsque l'on veut préciser $\OO$) de $\K^m$ tout $\OO$-module à droite engendré par une $\K$-base de $\K^m$. On l'appelle réseau \df{euclidien} si $\K=\R$, \df{hermitien complexe} si $\K=\C$ et \df{hermitien quaternionien} si $\K=\H$. On appelle \df{covolume} d'un réseau $\Gamma$ le volume du quotient $\K^m / \Gamma$ pour la mesure localement égale à la mesure de Haar sur $\K^m$, normalisée de sorte que le covolume du réseau standard $\OO^m$ soit égal à $1$. Un \df{réseau marqué} de $\K^m$ est un morphisme de $\OO$-modules à droite de $\OO^m$ dans $\K^m$ dont l'image est un réseau. Le groupe $\GL_m(\OO)$ agit à gauche par précomposition par l'adjoint, et le groupe de Lie $\GL_m(\K)$ agit à droite par postcomposition par l'adjoint sur l'ensemble des réseaux marqués de $\K^m$.

Notons $\mathcal{E}_m$ l'espace symétrique de type non compact du groupe de Lie quasi-simple $\SL_m(\K)$, c'est-à-dire l'espace homogène $\mathcal{E}_m = \SL_m(\K) / \SU_m(\K)$, où $\SU_m(\K)$ désigne le groupe spécial unitaire de $\K^m$, muni de l'action à gauche par translations de $\SL_m(\K)$ et d'une métrique riemannienne invariante par cette action.

Puisque $\OO$ est un ordre, tout morphisme de $\OO$-modules à droite de $\OO^m$ dans $\K^m$ s'étend uniquement en un endomorphisme de l'espace vectoriel à droite $\K^m$, et réciproquement tout endomorphisme de $\K^m$ se restreint en un morphisme de $\OO$-modules à droite de $\OO^m$ dans $\K^m$. Ainsi nous utiliserons les mêmes notations pour ces deux points de vue, sans risque de confusion.

Notons $\SL_m^1(\K)$ le sous-groupe de $\GL_m(\K)$ constitué des automorphismes de $\K^m$ dont le déterminant est de module $1$. L'ensemble des réseaux marqués de $\K^m$ de covolume $1$, muni de la topologie induite par la topologie produit sur $(\K^m)^{\OO^m}$, est muni de l'action à droite continue et simplement transitive de $\SL_m^1(\K)$ par postcomposition par l'adjoint, et est muni de l'action à gauche de $\SL_m(\OO)$ par précomposition par l'adjoint.

Notons $\mathcal{E}'_m$ l'ensemble des classes d'isométrie (positive ou non) de réseaux marqués de $\K^m$ de covolume $1$, muni de la topologie quotient : deux réseaux marqués $f$ et $f'$ sont identifiés s'il existe $g \in U_m(\K)$ tel que $g \circ f=f'$. L'homéomorphisme entre l'espace des réseaux marqués et $\SL_m^1(\K)$ passe au quotient en un homéomorphisme entre l'espace $\mathcal{E}'_m$ et l'espace homogène $\SL_m^1(\K) / U_m(\K) = \SL_m(\K) / \SU_m(\K) = \mathcal{E}_m$, qui est de plus $\SL_m(\OO)$-équivariant. Nous noterons dorénavant $\mathcal{E}_m$ tant l'espace des classes d'isométrie de réseaux marqués que l'espace symétrique, sans risque de confusion.

\bigskip

Rappelons ce qu'est une compactification d'un espace topologique $X$ localement compact : c'est la donnée d'une paire $(K,i)$, où $K$ est un espace topologique compact et $i:X \rightarrow K$ est un plongement d'image dense.  Si $G$ est un groupe agissant continûment sur $X$, on dit que $(K,i)$ est une \df{$G$-compactification} si l'action de $G$ sur $i(X)$, conjuguée par $i$ de l'action de $G$ sur $X$, s'étend continûment à $K$. Cette extension est alors unique. On dit que deux ($G$-) compactifications $(K,i)$ et $(K',i')$ de $X$ sont ($G$-) isomorphes s'il existe un homéomorphisme ($G$-équivariant) $f$ de $K$ sur $K'$ tel que $i'=f \circ i$.

\bigskip

Nous allons définir une compactification de Thurston des espaces symétriques $\mathcal{E}_m$, analogue à la compactification de Thurston des espaces de Teichmüller. Rappelons comment celle-ci est construite (voir par exemple \cite{thurston}, \cite{paulin_thurston}). Si $E$ est un ensemble, notons $\P({\R_+}^{\!E})$ l'espace topologique quotient de l'espace ${\R_+}^{\!E} \bs \{0\}$, muni de la topologie produit, par les homothéties de rapport strictement positif. Si $S$ est une surface compacte connexe orientée de genre supérieur ou égal à $2$ et si $\Gamma$ est son groupe fondamental, l'espace de Teichmüller de $S$ est l'ensemble des classes d'isométrie équivariante d'actions isométriques propres et libres de $\Gamma$ sur le plan hyperbolique réel $\H_\R^2$ (voir par exemple~\cite{paulin_bourbaki}). La compactification de Thurston de l'espace de Teichmüller de $S$ est alors l'adhérence de l'image du plongement dans $\P({\R_+}^{\!\Gamma})$, qui à la classe d'une telle action de $\Gamma$ sur $\H_\R^2$ associe la classe de l'application qui à un élément $\gamma$ de $\Gamma$ associe la distance de translation $\ell(\gamma)= \inf_{x \in \H_\R^2} d(x,\gamma x)$ de $\gamma$ dans $\H_\R^2$.

Nous allons modifier cette définition en remplaçant la surface $S$ par le tore $\K^m/\OO^m$, et le groupe $\Gamma$ par $\OO^m$, le groupe fondamental du tore. Notons $\phi$ l'application de $\mathcal{E}_m$ dans $\P({\R_+}^{\!\OO^m})$  qui à la classe d'isométrie équivariante d'un réseau marqué associe la classe d'homothétie de sa fonction distance de translation : un réseau marqué étant un morphisme de $\OO$-modules à droite $f$ de $\OO^m$ dans $\K^m$, notons $\phi([f])$ la classe d'homothétie de l'application
\beq \OO^m & \ra & \R_+ \\
u & \mapsto &  \| f(u) \| , \eeq
où $ \| . \| $ désigne la norme hermitienne de $\K^m$.

\blem L'application $\phi$ est bien définie et continue. \elem

Remarquons que puisque les espaces topologiques $\mathcal{E}_m$ et $\P({\R_+}^{\!\OO^m})$ sont métrisables, on peut utiliser les critères séquentiels pour montrer les propriétés topologiques de ces espaces.

\bp Composer un réseau marqué $f$ au but par une isométrie de $\K^m$ ne change pas $ \| f(u) \| $, pour tout $u \in \OO^m$. Donc l'application $\phi$ est bien définie.

Soit $(f_n)_{n \in \N}$ une suite de réseaux marqués de covolume $1$, qui converge vers un réseau marqué $f$ de covolume $1$. Alors la continuité de la norme assure que, pour tout $u \in \OO^m$, la suite $( \| f_n(u) \| )_{n \in \N}$ converge vers $ \| f(u) \| $. Ainsi la suite $(\phi([f_n]))_{n \in \N}$ converge vers $\phi([f])$. Donc l'application $\phi$ est continue. \ep

\blem \label{lem:image_phi} L'adhérence de l'image de $\phi$ est l'ensemble des classes d'homothétie d'applications
\beq \ell_f : \OO^m & \ra & \R_+ \\
u & \mapsto &  \| f(u) \| ,\eeq
où $f$ est un endomorphisme non nul de $\K^m$. L'image de $\phi$ est l'ensemble des classes d'homothétie d'applications $\ell_f$ pour lesquelles $f$ est inversible. \elem

\bp Soit $(f_n)_{n \in \N}$ une suite de réseaux marqués de covolume $1$, telle que la suite $(\phi([f_n]))_{n \in \N}$ converge vers la classe d'homothétie d'une application $\ell:\OO^m \ra \R_+$. 

Quitte à extraire, la suite $(f_n \K)_{n \in \N}$ converge vers $f \K$ dans l'espace projectif à droite des endomorphismes de $\K^m$, où $f$ est un endomorphisme de $\K^m$ non nul. Or, pour tout $v \in \K^m$ et tout $\lambda \in \K$, nous avons $\| v \lambda\| = \|v\| |\lambda|$. Alors, à homothétie réelle près, pour tout $u \in \OO^m$, la suite $( \| f_n(u) \| )_{n \in \N}$ converge vers $\ell(u)= \| f(u) \| $. Donc $\ell=\ell_f$ est bien du type décrit.

Réciproquement, soit $f$ un endomorphisme de $\K^m$ non nul, et soit $f_n = f+\frac{1}{n+1}\id$, pour tout entier $n \geq n_0$ tel que $f_n$ soit inversible. Alors l'endomoprhisme $|\det(f_n)|^\frac{1}{m}f_n$ appartient à $\SL_m^1(\K)$, donc définit un morphisme de $\OO$-modules de $\OO^m$ dans $\K^m$ dont l'image est de covolume $1$. Et la suite $(\phi([f_n]))_{n \geq n_0}$ converge vers la classe d'homothétie de l'application $\ell_f$ qui à $u \in \OO^m$ associe $ \| f(u) \| $ dans $\P({\R_+}^{\!\OO^m})$. 

\bigskip

Il est clair que l'image de $\phi$ est incluse dans l'ensemble des classes d'homothétie d'applications $\ell_f$ pour lesquelles $f$ est inversible. Réciproquement, si $f$ est un endomorphisme inversible de $\K^m$, alors à une homothétie réelle près on peut supposer que $|\det f|=1$, et donc $\R_+^* \ell_f = \phi([f])$, où $f : \OO^m \ra \K^m$ est bien un morphisme de $\OO$-modules dont l'image est de covolume $1$.
\ep

\blem L'image de $\phi$ est ouverte dans son adhérence. \elem

\bp Soit $f$ un morphisme de $\OO$-modules de $\OO^m$ dans $\K^m$ dont l'image est de covolume $1$. Soit $(f_n)_{n \in \N}$ une suite d'endomorphismes non nuls de $\K^m$, telle que la suite $(\R_+^*\ell_{f_n})_{n \in \N}$ converge vers $\phi([f])$. \`{A} des homothéties réelles près et quitte à extraire, on peut supposer que la suite $(f_n)_{n \in \N}$ converge vers un endomorphisme $g$ de $\K^m$ non nul. Dans ce cas, il est immédiat que la suite $(\ell_{f_n})_{n \in \N}$ converge vers $\ell_g$ dans ${\R_+}^{\!\OO^m}$. Donc $\R_+^* \ell_g = \phi([f])$, et puisque $f$ est inversible, la fonction $\ell_g$ ne s'annule qu'en $0 \in \OO^m$, et donc l'endomorphisme $g$ est inversible. L'ensemble des endomorphismes inversibles de $\K^m$ étant ouvert, on en déduit qu'à partir d'un certain rang les endomorphismes $f_n$ sont inversibles, et donc $\R_+^*\ell_{f_n}=\phi([f_n])$ appartient à l'image de $\phi$ : celle-ci est donc ouverte dans son adhérence.\ep

\blem L'adhérence de l'image de $\phi$ est compacte dans $\P({\R_+}^{\!\OO^m})$. \elem

\bp Soit $(\ell_n)_{n \in \N}$ une suite dans l'adhérence de l'image de $\phi$. D'après le lemme~\ref{lem:image_phi}, il existe pour tout $n \in \N$ un endomorphisme $f_n$ de $\K^m$ non nul tel que, pour tout $u \in \OO^m$, nous ayons $\ell_n(u)=\ell_{f_n}(u)= \| f_n(u) \| $. Par compacité de l'espace projectif à droite des endomorphismes de $\K^m$, quitte à extraire et à homothéties réelles près, on peut supposer que la suite $(f_n)_{n \in \N}$ converge vers un endomorphisme $f$ de $\K^m$ non nul. Alors la suite $(\R_+^*\ell_n)_{n \in \N}$ converge vers la classe d'homothétie de l'application $\ell_f$, et cette limite appartient à l'adhérence de l'image de $\phi$.
\ep

Pour montrer l'injectivité de l'application $\phi$, nous avons besoin d'un renforcement du lemme~\ref{lem:norme}.

\blem \label{lem:norme2} Soient $f$ et $f'$ deux endomorphismes de $\K^m$ tels que, pour tout $u \in \OO^m$, nous ayons $ \| f(u) \| = \| f'(u) \| $. Alors il existe $k \in \U_m(\K)$ tel que $kf=f'$. \elem

\bp Notons $d \in \{1,2,4\}$ la dimension de $\K$ sur $\R$, et soit $(e_1,\ldots,e_{dm})$ une $\R$-base de $\K^m$ formée par $m$ copies d'une $\Z$-base de $\OO$ : ainsi, tous les élements de cette base appartiennent à $\OO^m$. Soit $u \in \K^m$, et soit $u=\sum_{j=1}^{dm} e_j u_j$ sa décomposition dans cette base : les scalaires $u_j$ sont donc réels. Alors
$$ \| f(u) \| ^2 = \Big\| \sum_{j=1}^{dm} f(e_ju_j) \Big\| ^2 = \sum_{j,k=1}^{dm} \< f(e_j)u_j | f(e_k)u_k\> = \sum_{j,k=1}^{dm} u_j u_k \< f(e_j)|f(e_k)\>.$$
Or, par l'identité de polarisation pour le corps $\K$ (voir le lemme~\ref{lem:polarisation_C} ou \ref{lem:polarisation_H}), et par l'hypothèse sur $f$ et $f'$, nous savons que
$$ \forall j,k \in \lb 1,dm \rb, \;\; \< f(e_j)|f(e_k)\> \;\; = \;\; \< f'(e_j)|f'(e_k)\>. $$
On en déduit donc que $\|f(u)\| =  \| f'(u) \| $, et ceci pour tout $u \in \K^m$. D'après le lemme~\ref{lem:norme}, on en déduit qu'il existe $k \in \U_m(\K)$ tel que $kf=f'$. \ep

\bpro \label{pro:plongement} L'application $\phi$ est un plongement. \epro

\bp Montrons que l'application $\phi$ est injective : soient $f$ et $f'$ deux réseaux marqués de covolume $1$, tels que $\phi([f])=\phi([f'])$. Il existe donc un réel strictement positif $\lambda$ tel que pour tout $u \in \OO^m$, nous ayons $ \| f(u) \| = \lambda \| f'(u) \| = \| \lambda f'(u)\|$. D'après le lemme~\ref{lem:norme2}, on en déduit qu'il existe $k \in \U_m(\K)$ tel que $kf=\lambda f'$. Or $\det(kf) = \det(k) \det(f)$ est de module $1$, et $\det(\lambda f') = \lambda^m \det(f') = \lambda^m$ est un réel strictement positif, donc $\lambda = 1$. Ainsi $kf=f'$, donc $[f]=[f']$ dans $\mathcal{E}_m$ : l'application $\phi$ est injective.

\bigskip

Montrons que l'application $\phi$ est propre : soit $(f_n)_{n \in \N}$ une suite de réseaux marqués de covolume $1$, telle que la suite $(\phi([f_n]))_{n \in \N}$ converge vers $\phi([f])$, où $f$ est un réseau marqué de covolume $1$. Il existe donc une suite de réels strictement positifs $(\lambda_n)_{n \in \N}$ telle que, pour tout $u \in \OO^m$, la suite $(\lambda_n \|f_n(u)\|)_{n \in \N}$ converge vers $\|f(u)\|$. Ainsi la suite d'endomorphismes $(\lambda_n f_n)_{n \in \N}$ est bornée : quitte à extraire, on peut supposer qu'elle converge vers un endomorphisme $g$ de $\K^m$. Alors, pour tout $u \in \OO^m$, la suite $(\|\lambda_n f_n(u)\|)_{n \in \N}$ converge vers $\| g(u) \|  = \| f(u) \| $. D'après le lemme~\ref{lem:norme2}, on en déduit qu'il existe $k \in \U_m(\K)$ tel que $kg=f$, et en particulier $g$ a un déterminant de module $1$. Puisque la suite $(|\det(\lambda_n f_n)| = \lambda_n^m)_{n \in \N}$ converge vers $|\det g| = 1$, on en déduit que la suite $(\lambda_n)_{n \in \N}$ converge vers $1$. Donc la suite $(f_n)_{n \in \N}$ converge vers $g=k^{-1}f$, d'où la suite $([f_n])_{n \in \N}$ converge vers $[g] = [f]$ dans l'espace $\mathcal{E}_m$. L'application $\phi$ est ainsi propre.

\bigskip

L'application $\phi$, continue, injective et propre, est donc un plongement.
\ep

L'adhérence de l'image de $\phi$ dans $\P({\R_+}^{\!\OO^m})$ fournit une compactification de $\mathcal{E}_m$, donc de l'espace symétrique $\mathcal{E}_m$, que l'on appelle \df{compactification de Thurston}, et que l'on note $\ov{\mathcal{E}_m}^\mathcal{T}$. Notons que c'est une $\SL_m(\OO)$-compactifi\-cation, c'est-à-dire que l'action de $SL_m(\OO)$ s'étend contin\^ument au bord de $\mathcal{E}_m$, où $\SL_m(\OO)$ agit sur $\P({\R_+}^{\!\OO^m})$ à gauche par précomposition par l'adjoint et passage au quotient ${\R_+}^{\!\OO^m} \ra \P({\R_+}^{\!\OO^m})$.

\bigskip

Cette compactification est en fait munie d'une action continue à gauche de $\SL_m(\K)$ : pour le voir, on pourrait remplacer dans la construction qui précède $\OO^m$ par $\K^m$. Nous allons expliciter directement l'action de $\SL_m(\K)$ sur $\ov{\mathcal{E}_m}^\mathcal{T}$ : soit $g \in \SL_m(\K)$, et soit $\R_+^*\ell \in \ov{\mathcal{E}_m}^\mathcal{T}$. D'après le lemme~\ref{lem:image_phi}, il existe un endomorphisme $f$ de $\K^m$ non nul tel que, pour tout $u \in \OO^m$, nous ayons $\ell(u)= \| f(u) \| $. Définissons alors $g \cdot \R_+^*\ell$ comme la classe d'homothétie de l'application qui à $u \in \OO^m$ associe $ \| f(g^*(u)) \| $. C'est une action continue, qui étend l'action de $\SL_m(\OO)$ à gauche par précomposition par l'adjoint.

\section{Comparaison à la compactification de Satake}

Nous voulons comparer la compactification de Thurston définie ci-dessus à l'une des compactifications de Satake. Rappelons la construction de ces compactifications (voir~\cite{satake}). Soit $G$ un groupe de Lie réel connexe semi-simple de centre fini sans facteur compact, et soit $\rho : G \ra \SL(V)$ une représentation linéaire irréductible et de noyau fini de $G$ dans un $\K$-module à droite de dimension finie $V$. Soit $K$ un sous-groupe compact maximal de $G$ contenant $\Ker \rho$. Puisque $K$ est compact, il existe un produit scalaire hermitien sur $V$ tel que $\rho(K) \subset \SU(V)$. Si $h^*$ désigne l'adjoint d'un élément $h$ de $\SL(V)$ pour le produit scalaire hermitien de $V$, alors l'involution de Cartan de $\SL(V)$ associée à $\SU(V)$ est $h \mapsto (h^*)^{-1}$. Si $g \in G$ est tel que $\rho(g) = (\rho(g)^*)^{-1}$, alors $\rho(g) \in \SU(V)$, donc $g \in \rho^{-1}(SU(V)) = K$ car $K$ est maximal.

Notons $\P(\Sym(V))$ l'espace projectif de l'espace vectoriel réel des applications linéaires hermitiennes de $V$ dans $V$. Cet espace est muni d'une action naturelle de $G$ à gauche, définie par :
$$\forall g \in G, \forall \:\R h \in \P(\Sym(V)) \:,\: g \cdot \R h = \R \rho(g) h \rho(g)^*.$$

D'après~\cite{satake}, l'application de $G / K$ dans $\P(\Sym(V))$ qui à $gK$ associe $\R \rho(g)\rho(g)^*$ est un plongement, dont l'adhérence de l'image est appelée la \df{compactification de Satake} de $G/K$ associée à la représentation $\rho$. L'action de $G$ sur $\P(\Sym(V))$ préserve l'image de ce plongement, qui est de plus équivariant pour les actions de $G$ : c'est donc une $G$-compactification de $G/K$.

\bigskip

Considérons ici le groupe $G=\SL_m(\K)$, et le sous-groupe $K=\SU_m(\K)$. Considérons la représentation $\rho = \id$ de $\SL_m(\K)$ pour l'action linéaire à gauche de $\SL_m(\K)$ sur le $\K$-module à droite $V=\K^m$. Considérons le plongement associé
\beq  \SL_m(\K) / \SU_m(\K)& \ra & \P(\Sym(\K^m)) \\
\gamma \SU_m(\K) & \mapsto & \R \gamma \gamma^*. \eeq
Notons $\ov{\mathcal{E}_m}^\mathcal{S}$ l'adhérence de son image, c'est-à-dire la compactification de Satake de $\mathcal{E}_m$ associée à la représentation $\rho$.

\bthm \label{thm:thurston_satake} Les deux compactifications $\ov{\mathcal{E}_m}^\mathcal{T}$ et $\ov{\mathcal{E}_m}^\mathcal{S}$ sont $\SL_m(\K)$-iso\-morphes. \ethm

\bp La compactification de Satake $\ov{\mathcal{E}_m}^\mathcal{S}$, adhérence des classes d'homothétie des matrices hermitiennes définies positives, est l'ensemble des classes d'homothétie des matrices hermitiennes positives non nulles : en effet, soit $a$ une matrice hermitienne positive non nulle. Alors, pour $n \geq n_0$, la matrice $a+\frac{1}{n+1}I_m$ est hermitienne, positive et inversible, donc définie positive. Ainsi $(a+\frac{1}{n+1}I_m)_{n \geq n_0}$ est une suite de matrices hermitiennes définies positives convergeant vers $a$. Par ailleurs, l'ensemble des matrices hermitiennes positives est fermé. Nous identifierons donc $\ov{\mathcal{E}_m}^\mathcal{S}$ avec l'espace des classes d'homothétie de matrices hermitiennes positives non nulles.

Définissons une application $\xi : \ov{\mathcal{E}_m}^\mathcal{S} \ra \P({\R_+}^{\!\OO^m})$. Soit $a$ une matrice hermitienne positive non nulle. Par diagonalisation (voir le lemme~\ref{lem:diag} si $\K=\H$), la matrice $a$ admet une unique racine carrée hermitienne positive que l'on notera $\sqrt{a}$, c'est-à-dire telle que $\sqrt{a}^2=a$. Remarquons que si $\lambda \in \R_+$, alors $\sqrt{\lambda a} = \sqrt{\lambda} \sqrt{a}$. Définissons l'application
\beq \widehat{\xi}(a) : \OO^m & \ra & \R_+ \\
u & \mapsto &  \| \sqrt{a}(u) \|, \eeq
ce qui permet de poser
\beq \xi : \ov{\mathcal{E}_m}^\mathcal{S} &\ra& \P({\R_+}^{\!\OO^m}) \\
\R a &\mapsto& \R_+^* \widehat{\xi}(a),\eeq
application qui est bien définie.

\bigskip

Montrons que l'image de $\xi$ est égale à $\ov{\mathcal{E}_m}^\mathcal{T}$. D'après la définition de $\xi$ et le lemme~\ref{lem:image_phi}, elle est incluse dedans. Maintenant, soit $f$ un endomorphisme non nul de $\K^m$, et soit $\ell_f: u \mapsto  \| fu \| $. Montrons alors que $\widehat{\xi}(f^*f)=\ell_f$, ce qui conclura par le lemme~\ref{lem:image_phi}. Soit $u \in \OO^m$, et soit $f=kp$ une décomposition polaire de $f$ : $k \in \U_m(\K)$ et $p$ est un endomorphisme hermitien positif non nul. Alors $\sqrt{f^*f}=\sqrt{p^*k^*kp}=\sqrt{p^2}=p$, donc $\widehat{\xi}(f^*f)(u)= \| \sqrt{f^*f}u \| = \| pu \| = \| kpu \| = \| fu \| = \ell_f(u)$.

\bigskip

Cette application est $\SL_m(\K)$-équivariante : soit $a$ une matrice hermitienne positive non nulle, et soit $g \in \SL_m(\K)$. D'après la décomposition polaire (voir le lemme~\ref{lem:polaire} si $\K=\H$), il existe une matrice unitaire $k$ et une matrice hermitienne positive $p$ non nulle telles que $g\sqrt{a}=pk$. Soit $u \in \OO^m$, alors 
$$ \widehat{\xi}(g \cdot a)(u) = \widehat{\xi}(gag^*)(u) =  \|  \sqrt{gag^*}(u)  \| .$$
Or $gag^*= g\sqrt{a} (g\sqrt{a})^* = (pk)(pk)^*= p p^* = p^2$, donc $\sqrt{gag^*}=p$. Ainsi
$$ \widehat{\xi}(g \cdot a)(u) =  \| pu \|  =  \| p^* u \|  =  \| k\sqrt{a}g^*u \|  =  \| \sqrt{a}g^*u \|  = (g \cdot \widehat{\xi}(a))(u),$$
où nous renvoyons à la fin de la partie~\ref{sec:thurston} pour la définition de l'action de $\SL_m(\K)$ sur $\ov{\mathcal{E}_m}^\mathcal{T}$. Nous avons ainsi $\xi(g \cdot \R a) = g \cdot \xi(\R a)$.

\bigskip

L'application qui à une matrice hermitienne positive $a$ associe sa racine carrée hermitienne positive $\sqrt{a}$ est continue, et par continuité de la norme on en déduit que, à $u \in \OO^m$ fixé, l'application qui à $a$ associe $\widehat{\xi}(a)(u)$ est continue. Puisque l'espace $\P({\R_+}^{\!\OO^m})$ est muni de la topologie quotient de la topologie produit sur ${\R_+}^{\!\OO^m}$, il s'ensuit que l'application $\widehat{\xi}$ est continue. On en déduit que l'application quotient $\xi$ est continue.

\bigskip

Montrons enfin que l'application $\xi$ est injective : soient $a$ et $a'$ deux matrices hermitiennes positives non nulles telle que $\xi(\R a)=\xi(\R a')$. \`{A} une homothétie de rapport strictement positif près, on peut supposer que $\widehat{\xi}(a)=\widehat{\xi}(a')$. Ainsi, pour tout $u \in \OO^m$, nous avons $ \| \sqrt{a}u \| = \| \sqrt{a'}u \| $. D'après le lemme~\ref{lem:norme2}, on en déduit qu'il existe $k \in \U_m(\K)$ tel que $k \sqrt{a} = \sqrt{a'}$. Par unicité dans la décomposition polaire, on conclut que $\sqrt{a}=\sqrt{a'}$, d'où $a=a'$. Ainsi l'application $\xi$ est injective.

\bigskip

Ainsi l'application $\xi$ est une bijection continue et $\SL_m(\K)$-équivariante de l'espace compact $\ov{\mathcal{E}_m}^\mathcal{S}$ sur l'espace séparé $\ov{\mathcal{E}_m}^\mathcal{T}$, donc est un homéomorphisme $\SL_m(\K)$-équivariant.
\ep

\section{Compactification d'espaces de réseaux autoduaux}

\label{sec:autoduaux}

Pour les prérequis de cette partie, on renvoie à~\cite{bavard}. Fixons $\tau = \id_\R$ si $\K=\R$, $\tau = \id_\C$ ou la conjugaison si $\K=\C$, et $\tau$ la conjugaison si $\K=\H$. Fixons $b$ une forme $\tau$-sesquilinéaire à gauche sur $\K^m$ non dégénérée, hermitienne ou anti-hermitienne. Dans une base adaptée, la forme $b$ est définie par
$$ \forall x,y \in \K^m \:,\: b(x,y) = \,\t x^\tau Jy,$$
où $J \in U_m(\K) \cap \GL_m(\OO)$.

Notons $\SU(b)$ le sous-groupe de $\GL_m(\K)$ constitué des automorphismes de $\K^m$ qui préservent $b$ : puisque $J$ est unitaire, le groupe $\SU(b)$ est autoadjoint.  Si $\Lambda$ est un $\OO$-réseau de $\K^m$, on définit son dual par rapport à $b$ par :
$$ \Lambda^{*b} = \{y \in \K^m \,:\, \forall x \in \Lambda, \, b(x,y) \in \OO\}.$$
On dit que le réseau $\Lambda$ est \df{autodual} (pour $b$) si $\Lambda^{*b}=\Lambda$. Par exemple, le réseau $\OO^m$ est autodual.
%
%
%
%
%
%
%
%

Soit $\Lambda_0 = \OO^m$ le $\OO$-réseau standard autodual de covolume $1$, marqué par l'identité $f_0 : \OO^m \ra \Lambda_0$.

Pour l'action à droite de $\SL_m(\K)$ sur l'espace des réseaux marqués par postcomposition par l'adjoint, l'orbite du réseau marqué $f_0$ par $\SU(b)$ est constituée de réseaux autoduaux. En effet, soit $g \in \SU(b)$, et considérons $y \in (g^* \Lambda_0)^{*b}$. Nous savons donc que $\forall x \in g^* \Lambda_0, b(x,y) \in \OO$. Ainsi $\forall x \in \Lambda_0, b(g^*x,y) = b(x,{g^*}^{-1}y) \in \OO$. Donc ${g^*}^{-1}y \in \Lambda_0^{*b} = \Lambda_0$. Ainsi $(g^* \Lambda_0)^{*b} = g^* \Lambda_0$. 

Dans l'identification entre l'espace $\mathcal{E}_m$ des classes d'isométrie de réseaux marqués de covolume $1$ et l'espace symétrique $\SL_m(\K) / \SU_m(\K)$, l'ensemble des classes d'isométrie de réseaux marqués dans l'orbite de $f_0$ sous $\SU(b)$ s'identifie alors à l'espace homogène $\SU(b) \SU_m(\K) / \SU_m(\K)$, c'est-à-dire à l'espace $\SU(b) / (\SU(b) \cap \SU_m(\K))$. Notons $\mathcal{E}^b_m$ ce sous-espace de $\mathcal{E}_m$.

On supposera que le groupe $\SU(b)$ n'est pas compact, de sorte que l'espace homogène $\mathcal{E}^b_m$ n'est pas compact. Lorsque l'on restreint le plongement $\phi$ de $\mathcal{E}_m$ dans ${\mathcal{E}_m}^\mathcal{T}$ au sous-espace $\mathcal{E}^b_m$, l'adhérence de son image définit la \df{compactification de Thurston} $\ov{\mathcal{E}^b_m}^\mathcal{T}$ de l'espace $\mathcal{E}^b_m$.

Lorsque l'on restreint le plongement $\mathcal{E}_m$ dans ${\mathcal{E}_m}^\mathcal{S}$ au sous-espace $\mathcal{E}^b_m$, l'adhérence de son image définit la \df{compactification de Satake} $\ov{\mathcal{E}^b_m}^\mathcal{S}$ de l'espace $\mathcal{E}^b_m$. 

Le groupe $\SU(b)$ est un groupe de Lie réel simple non compact. Alors la compactification de Satake définie ci-dessus est naturellement isomorphe à la compactification de Satake associée à la représentation $\rho$ de $\SU(b)$ dans $\SL_m(\K)$, restriction de l'identité de $\SL_m(\K)$.

\bpro La compactification de Thurston $\ov{\mathcal{E}^b_m}^\mathcal{T}$ de l'espace des réseaux marqués autoduaux est $\SU(b)$-isomorphe à la compactification de Satake $\ov{\mathcal{E}^b_m}^\mathcal{S}$. \epro

\bp L'homéomorphisme $\SL_m(\K)$-équivariant (voir le théorè\-me~\ref{thm:thurston_satake}) entre $\ov{\mathcal{E}_m}^\mathcal{T}$ et $\ov{\mathcal{E}_m}^\mathcal{S}$ induit un homéomorphisme $\SU(b)$-équivariant entre $\ov{\mathcal{E}^b_m}^\mathcal{T}$ et $\ov{\mathcal{E}^b_m}^\mathcal{S}$.
\ep

Dans la suite, nous aurons besoin de savoir que dans le cas de l'espace vectoriel réel $\R^{2g}$ muni de la forme symplectique standard
$$ b(x,x') = \sum_{j=1}^g x_jx'_{g+j} - x_{g+j}x'_{j},$$
tout réseau autodual pour $b$ est dans l'orbite du réseau $\Z^{2g}$.

\blem \label{lem:type_symplectique} Tout réseau autodual pour $b$ de $\R^{2g}$ est l'image par un élément de $\Sp_{2g}(\R)$ du réseau standard $\Z^{2g}$. \elem

\bp Soit $\Lambda$ un réseau autodual de $\R^{2g}$ de covolume $1$. Soit $A \in \SL_{2g}(\R)$ une matrice telle que $\Lambda = A \cdot \Z^{2g}$. Le réseau $\Lambda$ étant autodual pour la forme symplectique $b$ qui a pour matrice
$$ J = \left( \begin{array}{cc} 0 & -I_g \\ I_g & 0 \end{array} \right),$$
on en déduit que le réseau $\t AJA \cdot \Z^{2g}$ est dual du réseau $\Z^{2g}$ pour le produit scalaire standard sur $\R^{2g}$. Or on sait que le réseau $\Z^{2g}$ est autodual pour le produit scalaire standard, donc $\t AJA \cdot \Z^{2g} = \Z^{2g}$ : en particulier $\t AJA \in SL_{2g}(\Z)$. C'est la matrice d'une forme symplectique à coefficients entiers, donc par réduction il existe une matrice $B \in \GL_{2g}(\Q)$ telle que $\t B\t AJAB = J$. Ainsi la matrice $AB$ appartient à $\Sp_{2g}(\R)$.

Quitte à prendre l'image du réseau $\Lambda$ par $(AB)^{-1} \in \Sp_{2g}(\R)$, on suppose que le réseau $\Lambda = (AB)^{-1} A \cdot \Z^{2g} = B^{-1} \cdot \Z^{2g}$ est commensurable à $\Z^{2g}$. Le fait que $\Lambda$ soit autodual impose alors qu'il s'écrive
$$ \Lambda = \bigoplus_{j=1}^g \left( \Z r_j e_j \oplus \Z \frac{1}{r_j} e_{g+j}\right),$$
où $r_1, \ldots,r_g$ sont des rationnels non nuls. Puisque la matrice
$$C=\Diag(r_1^{-1},\ldots,r_g^{-1},r_1,\ldots,r_g)$$
appartient à $\Sp_{2g}(\R)$ et que $C \cdot \Lambda = \Z^{2g}$, ceci conclut le lemme. \ep

\section{Cas du groupe de Lie exceptionnel $E_{6(-26)}$}

Notons $\O$ désigne l'algèbre non associative des octonions de Cayley (voir \cite{baez,allcock_cayley,mare_willems,octonions}). Le groupe $\SL_3(\O)$ (dont nous rappelons la définition ci-dessous) est une forme réelle non compacte de rang réel $2$ du groupe de Lie complexe exceptionnel $E_6$, notée $E_{6(-26)}$. Montrons comment les résultats précédents s'étendent à ce groupe.

L'algèbre non associative $\O$ des octonions de Cayley est l'espace vectoriel réel euclidien de dimension $8$, de base orthonormée $(e_0=1,e_1,\ldots,e_7)$, muni de la multiplication bilinéaire définie sur cette base par la table~\ref{tab:octonions}.

\begin{table}[!h]
$$\begin{array}{|c|c|c|c|c|c|c|c|}
\hline
& e_1 & e_2 & e_3 & e_4 & e_5 & e_6 & e_7  \\
\hline
e_1 & -1 & e_4 & e_7 & -e_2 & e_6 & -e_5 & -e_3  \\
\hline
e_2 & -e_4 & -1 & e_5 & e_1 & -e_3 & e_7 & -e_6  \\
\hline
e_3 & -e_7 & -e_5 & -1 & e_6 & e_2 & -e_4 & e_1  \\
\hline
e_4 & e_2 & -e_1 & -e_6 & -1 & e_7 & e_3 & -e_5  \\
\hline
e_5 & -e_6 & e_3 & -e_2 & -e_7 & -1 & e_1 & e_4  \\
\hline
e_6 & e_5 & -e_7 & e_4 & -e_3 & -e_1 & -1 & e_2  \\
\hline
e_7 & e_3 & e_6 & -e_1 & e_5 & -e_4 & -e_2 & -1  \\
\hline
\end{array}
$$
\caption{Table de multiplication des octonions}
\label{tab:octonions}
\end{table}

C'est une algèbre non associative à division, munie de la conjugaison
\beq \O & \ra & \O \\
x = x_0e_0 + \sum_{i=1}^7 x_i e_i & \mapsto & \ov{x} = x_0e_0 - \sum_{i=1}^7 x_i e_i.\eeq
La norme euclidienne de $\O$ vérifie $\|x\| = \sqrt{\ov{x}x} = \sqrt{x\ov{x}}$ et $\|xy\|=\|x\| \|y\|$. Pour tout octonion $x \in \O$, on définit de plus sa partie réelle $\Re x = \frac{x+\ov{x}}{2}$ et sa partie imaginaire $\Im x = \frac{x-\ov{x}}{2}$.

Soit $m$ un entier au moins égal à $2$. L'espace vectoriel réel $\O^m$ est naturellement muni d'une structure euclidienne, pour le produit scalaire
$$ \<u,v\> = \sum_{i=1}^m \Re(\ov{u_i}v_i) \;,\; \mbox{ où } u,v \in \O^m.$$
Considérons $\M_m(\O)$ l'espace vectoriel réel des matrices carrées de taille $m$ à coefficients dans $\O$. En tant qu'ensemble d'endomorphismes $\R$-linéaires de $\O^m$, muni de la composition des endomorphismes, c'est une algèbre associative. Par contre, la multiplication n'est pas obtenue avec la formule usuelle du produit matriciel. Considérons le groupe $\GL_m(\O)$ des matrices de $\M_m(\O)$ qui induisent un isomorphisme de $\O^m$.

L'adjoint $M^*$ d'une matrice $M \in \M_m(\O)$ pour le produit scalaire de $\O^m$ est alors la matrice transposée et conjuguée de $M$. On dit que la matrice $M$ est hermitienne si $M=M^*$, et on note $\hh_m(\O)$ l'espace vectoriel réel des matrices hermitiennes. Une matrice hermitienne $M$ est dite positive (resp. définie positive) si pour tout $x \in \O^m \bs \{0\}$, nous avons $\<x,Mx\> \geq 0$ (resp. $\<x,Mx\> > 0$).

Une matrice $M \in \M_m(\O)$ est dite unitaire si $MM^*=I_m$. Le sous-groupe $U_m(\O)$ de $\GL_m(\O)$ constitué des matrices unitaires est un sous-groupe compact maximal de $\GL_m(\O)$.

Les résultats d'algèbre linéaire de la partie~\ref{sec:algebre_lineaire} se généralisent dans ce cadre, montrons comment.

\bigskip

Le sous-espace vectoriel de $\O$ engendré par $1,e_1,e_2$ et $e_4$ est une algèbre associative isomorphe au corps gauche $\H$ : nous identifierons ainsi $i$ avec $e_1$, $j$ avec $e_2$ et $k$ avec $e_4$. Une base réelle de $\O$ est alors donnée par $(1,i,j,k,e_3,ie_3,je_3,ke_3)$. Le sous-espace vectoriel de $\O$ engendré par $1$ et $e_3$ est une algèbre commutative isomorphe à $\C$, et nous identifierons les deux : une $\C$-base de $\O$ est alors donnée par $(1,i,j,k)$. Considérons l'isomorphisme $\C$-linéaire à droite
\beq \beta : \O^m & \ra & \C^{4m} = (\C^4)^m \\
(x_l + iy_l + jz_l + kw_l)_{l \in \lb 1,m \rb} & \mapsto & (x_l,y_l,z_l,w_l)_{l \in \lb 1,m \rb}.\eeq
Considérons alors l'application
\beq \eta : \M_m(\O) & \ra & \M_{4m}(\C) \\
M & \mapsto & \left( u \in \C^{4m} \mapsto \beta(M(\beta^{-1}(u))) \in \C^{4m}\right),\eeq
c'est un plongement de $\C$-algèbres à droite associatives, équivariant pour l'adjoint. L'application $\eta$ réalise le sous-groupe $\GL_m(\O)$ comme un sous-groupe algébrique autoadjoint de $\GL_{4m}(\C)$ : d'après~\cite[Theorem~7.1, p.~224]{helgason}, l'orbite de ce sous-groupe dans l'espace symétrique de $\GL_{4m}(\C)$ est un sous-espace symétrique, donc la décomposition polaire et la décomposition de Cartan y sont internes. En particulier, la diagonalisation des matrices hermitiennes inversibles y est interne, donc par densité on en déduit la diagonalisation des matrices hermitiennes quelconques. Nous pouvons donc énoncer les deux résultats suivants.

\blem[Diagonalisation des matrices hermitiennes octonioniques] \label{lem:diag_O} Pour toute matrice hermitienne $M \in \hh_m(\O)$, il existe une matrice unitaire $U \in U_m(\O)$ et une matrice diagonale réelle $D \in M_m(\R)$ telles que $M=UDU^{-1}$.
\elem

\blem[Décomposition polaire octonionique] \label{lem:polaire_O} Soit $M \in \M_m(\O)$. Alors il existe une matrice unitaire $U \in U_m(\O)$ et une unique matrice hermitienne positive $P \in \hh_m(\O)$ telles que $M = PU$. \elem

On dispose également d'une identité de polarisation.

\blem[Identité de polarisation octonionique] \label{lem:polarisation_O} Soit $(q_1,\ldots,q_8)$ une $\R$-base de $\O$. Alors il existe huit réels $(\lambda_1,\ldots,\lambda_8) \in \R^8$ tels que, pour tous $u,v \in \O^m$, nous ayons :
$$ \<u|v\> = \sum_{l=1}^8 \left(  \| u+ vq_l \| ^2 -  \| u- vq_l \| ^2\right) \lambda_l .$$
\elem

\bp Remarquons tout d'abord que
$$ \| u+ v q_l \| ^2 -  \| u- v q_l \| ^2 = 4 \Re(\ov{u}(vq_l)).$$
Or la partie réelle des octonions est associative, donc $\| u+ v q_l \| ^2 -  \| u- v q_l \| ^2 = 4 \Re((\ov{u}v)q_l)$.
Considérons l'application $\R$-linéaire
\beq \varphi : \O & \ra & \R^8 \\
q & \mapsto & (4\Re(qq_l))_{l \in \lb 1,8 \rb}.\eeq
Si $q \in \O$ appartient au noyau de $\varphi$, alors par $\R$-linéarité $4\Re(q\ov{q}) = 4|q|^2 = 0$ donc $q=0$ : l'application $\varphi$ est donc un isomorphisme $\R$-linéaire. Pour tout $l \in \lb 1,8 \rb$, notons $\lambda_l = \Re(\varphi^{-1}(f_l))$, où $(f_1,\ldots,f_8)$ désigne la base canonique de $\R^8$.
Alors, pour tous $u,v \in \O^m$, nous avons par $\R$-linéarité
\beq \<u,v\> &=& \Re(\ov{u}v) = \Re \left(\varphi^{-1}(\varphi(\ov{u}v))\right) \\
&=& \Re \left(\varphi^{-1}\left((4\Re((\ov{u}v) q_l))_{l \in \lb 1,8 \rb}\right)\right) \\
&=& \Re \left(\varphi^{-1}\left( \sum_{l=1}^8 4\Re((\ov{u}v)q_l) f_l \right)\right) \\
&=& \sum_{l=1}^8 4\Re((\ov{u}v)q_l) \Re(\varphi^{-1}(f_l)) \\
\hspace{2.7cm} &=& \sum_{l=1}^8 \left(  \| u+ v q_l \| ^2 -  \| u- v q_l \| ^2\right) \lambda_l. \hspace{2.7cm} \qed \eeq
\epp

Définissons $\SL_m(\O)$ le sous-groupe de $\GL_m(\O)$ constitué des matrices qui préservent une mesure de Haar sur $\O^m$. On peut également définir le déterminant d'une matrice hermitienne, auquel cas c'est aussi le sous-groupe qui préserve le déterminant, pour l'action suivante :
$$ \forall g \in \GL_m(\O), \forall M \in \hh_m(\O),\; g \cdot M = gM+Mg^*.$$
Voici les formules définissant le déterminant des matrices hermitiennes dans les cas $m=2$ et $m=3$.
\bit
\item Lorsque $m=2$, les matrices hermitiennes s'écrivent
$$ \hh_2(\O) = \left\{ \left( \begin{array}{cc} \alpha & x \\ \ov{x} & \beta \end{array} \right) \;,\; x \in \O, \alpha, \beta \in \R \right\}.$$
On définit alors le déterminant par la formule
$$ \det \left( \begin{array}{cc} \alpha & x \\ \ov{x} & \beta \end{array} \right) = \alpha \beta - |x|^2.$$
\item Lorsque $m=3$, les matrices hermitiennes s'écrivent
$$ \hh_3(\O) = \left\{ \left( \begin{array}{ccc} \alpha & \ov{z} & \ov{y} \\ z & \beta & x \\ y & \ov{x} & \gamma \end{array} \right) \;,\; x,y,z \in \O, \alpha, \beta,\gamma \in \R \right\}.$$
On définit alors le déterminant par la formule
$$ \det \left( \begin{array}{ccc} \alpha & \ov{z} & \ov{y} \\ z & \beta & x \\ y & \ov{x} & \gamma \end{array} \right) = \alpha \beta \gamma - \left(\alpha|x|^2+\beta|y|^2 + \gamma|z|^2\right) +2\Re(xyz).$$
\eit

Lorsque $m=2$, le groupe de Lie réel $\SL_2(\O)$ est isomorphe au groupe simple $\Spin(9,1)$, revêtement universel du groupe $\SO_0(9,1)$. C'est un groupe classique, dont l'espace symétrique est l'espace des classes d'isométrie de réseaux orthogonaux de $\R^{10}$ pour une forme quadratique de signature $(9,1)$, donc son cas est traité dans la partie~\ref{sec:autoduaux}.

Dorénavant nous nous placerons dans le cas $m=3$. Alors le groupe de Lie simple $G=\SL_3(\O)$ est une forme réelle non compacte de rang réel $2$ du groupe de Lie complexe exceptionnel $E_6$, notée $E_{6(-26)}$. Le sous-groupe $K = \SU_3(\O) = \SL_3(\O) \cap U_3(\O)$ est un sous-groupe compact maximal de $\SL_3(\O)$, isomorphe au groupe $F_{4(-52)}$, la forme réelle compacte du groupe de Lie complexe exceptionnel $F_4$.

Notons $\OO$ un ordre (non associatif) de $\O$, c'est-à-dire le sous-groupe additif engendré par une base réelle de $\O$ et stable par multiplication, contenant l'anneau des entiers $\Z[e_0,\ldots,e_7]$ (il existe $16$ tels ordres, voir par exemple~\cite[Theorem1, p.~100]{octonions}). Appelons $\OO$-réseau marqué de $\O^3$ toute application $\OO$-équivariante à droite $f: \OO^3 \ra \O^3$ telle que le $\R$-espace vectoriel engendré par l'image soit égal à $\O^3$. Le covolume d'un réseau marqué est le covolume de son image dans $\O^3$, pour une mesure de Haar sur $\O^3$ normalisée de sorte que le covolume de $\OO^3$ soit égal à $1$.

\`{A} toute matrice hermitienne $M$ de $\hh_3(\O)$ définie positive de déterminant $1$, on associe la classe d'isométrie du $\OO$-réseau marqué $f :\OO^3 \ra \O^3$ de covolume $1$ qui à $u \in \OO^3$ associe $\| M^*f(u)\|$. Et à toute classe d'isométrie d'un $\OO$-réseau marqué $f :\OO^3 \ra \O^3$ de covolume $1$, on associe la classe à droite modulo $\SU_3(\O)$ d'une matrice $M \in \SL_3(\O)$ telle que pour tout $u \in \OO^3$ nous ayons $f(u) = M^*(u)$, donc d'après la décomposition polaire de $G$ on associe une unique matrice hermitienne de $\hh_3(\O)$ définie positive de déterminant $1$. Ainsi l'espace $X$ des classes d'isométrie de $\OO$-réseaux marqués de covolume $1$ de $\O^3$ s'identifie au sous-espace de $\hh_3(\O)$ constitué des matrices hermitiennes définies positives de déterminant $1$, et ce sont deux modèles de l'espace symétrique du groupe $G$.

On peut alors considérer l'application
\beq \phi : X & \mapsto & \P({\R_+}^{\!\OO^3}) \\
\ [f:\OO^3 \ra \O^3] & \ra & [u \mapsto \|f(u)\|]. \eeq

On peut alors appliquer la même preuve que celle de la proposition~\ref{pro:plongement}, en remplaçant \og endomorphisme \fg\ par \og endomorphisme $\R$-linéaire $\OO^3$-équi\-variant à droite \fg. On en déduit que l'application $\phi$ est un plongement, et on appelle l'adhérence de son image la compactification de Thurston $\ov{X}^\mathcal{T}$ de $X$.

Par ailleurs, la compactification de Satake de $X$ associée à la représentation $\rho : G \ra \End_\R(\hh_3(\O))$ définie par $g \cdot M = gM + Mg^*$ est l'adhérence $\ov{X}^\mathcal{S}$ de l'image de $X$ dans l'espace projectif réel $\P(\hh_3(\O))$ de $\hh_3(\O)$.

La preuve du théorème~\ref{thm:thurston_satake} est vraie dans ce cadre, et permet d'énoncer le résultat suivant.

\bthm \label{thm:thurston_satake_O} Les deux compactifications $\ov{X}^\mathcal{T}$ et $\ov{X}^\mathcal{S}$ de l'espace symétrique $X$ du groupe de Lie exceptionnel $E_{6(-26)}$ sont $\SL_3(\O)$-isomorphes. \ethm

\section{Compactification de Thurston de l'espace de Torelli}

\label{sec:torelli}

Le but de cette partie est de construire une compactification de l'espace de Torelli d'une surface $S$ analogue à la compactification de Thurston de l'espace de Teichmüller de $S$. Pour les définitions de base, nous renvoyons à~\cite{mcg}, \cite{bost} et \cite{warner_diff}.

Soient $g>0$ et $q \geq 0$ deux entiers tels que $2-2g-q < 0$. Fixons $S$ une surface lisse compacte connexe orientée de genre $g$, munie d'une partie fixée $P$ de cardinal $q$, dont les éléments sont appelés les points marqués de $S$. Dans le cas où $q=0$, on omettra la notation $P$ dans ce qui suit.

\bigskip

Notons $\Diff^+(S,P)$ le groupe des difféomorphismes de $S$ préservant l'o\-rientation, fixant $P$ point par point, muni de la topologie compacte-ouverte.

Notons $\Diff_{H^1}(S,P)$ le sous-groupe des difféomorphismes de $S$ préservant l'orientation, fixant $P$ point par point, induisant l'identité sur $H^1(S \bs P,\R)$.

Notons $\Diff_0(S,P)$ le sous-groupe des difféomorphismes de $S$ préservant l'orientation, isotopes à l'identité par une isotopie fixant $P$.

Nous avons les inclusions $\Diff_0(S,P) \subset \Diff_{H^1}(S,P) \subset \Diff^+(S,P)$.

\bigskip

Notons $\Teich(S,P)$ l'espace de Teichmüller de $S$ : c'est l'ensemble des classes d'équivalence de couples $(X,h)$, où $X$ est une surface hyperbolique orientée complète d'aire finie et $h: S \bs P \ra X$ est un difféomorphisme préservant l'orientation (appelé un marquage), et où l'on identifie les couples $(X,h)$ et $(X',h')$ s'il existe une isométrie préservant l'orientation $s:X \ra X'$ telle que $h'$ soit égal à $s \circ h$ modulo $\Diff_0(S,P)$, c'est-à-dire que telle que $h'$ soit isotope à $s \circ h$. Le groupe $\Diff^+(S,P)$ agit sur $\Teich(S,P)$ par précomposition du marquage, de noyau d'action $\Diff_0(S,P)$. Le groupe modulaire $\MCG(S,P) = \Diff^+(S,P) / \Diff_0(S,P)$ agit ainsi fidèlement sur $\Teich(S,P)$.

Notons $\Tor(S,P)$ l'espace de Torelli de $S$ : c'est l'ensemble des classes d'équivalence de couples $(X,h)$, où $X$ est une surface hyperbolique orientée complète d'aire finie et $h: S \bs P \ra X$ est un difféomorphisme préservant l'orientation, et où l'on identifie les couples $(X,h)$ et $(X',h')$ s'il existe une isométrie préservant l'orientation $s:X \ra X'$ telle que $h'$ soit égal à $s \circ h$ modulo $\Diff_{H^1}(S,P)$. Ceci revient à demander l'égalité $h'^* = h^* \circ s^* : H^1(X',\R) \ra H^1(S \bs P,\R)$.

Notons $\Mod(S,P)$ l'espace des modules de $S$ : c'est l'ensemble des classes d'équivalence de couples $(X,h)$, où $X$ est une surface hyperbolique  orientée complète d'aire finie et $h: S \bs P \ra X$ est un difféomorphisme préservant l'orientation, et où l'on identifie les couples $(X,h)$ et $(X',h')$ si $X$ et $X'$ sont isométriques.

\bigskip

Munissons l'espace de Teichmüller $\Teich(S,P)$ de la topologie quotient de la topologie induite par celle de $\Gamma(\otimes^2 TS)$, où à chaque couple $(X,h)$ on associe $h^* \sigma_X$, en notant $\sigma_X$ la métrique riemannienne hyperbolique de $X$. L'espace de Teichmüller admet alors une structure naturelle de variété complexe biholomorphe à $\C^{3g-3+q}$ (voir~\cite{thurston}).

L'espace de Torelli $\Tor(S,P)$ est le quotient de $\Teich(S,P)$ par le groupe de Torelli $T(S,P) = \Diff_*^+(S,P) / \Diff_0^+(S,P)$. Le groupe de Torelli étant sans torsion (voir par exemple~\cite{mcg}), on munit $\Tor(S,P)$ de la structure de variété quotient induite par celle de $\Teich(S,P)$.

L'espace des modules $\Mod(S,P)$ est le quotient de $\Teich(S,P)$ par le groupe modulaire $\MCG(S,P)$. On munit $\Mod(S,P)$ de la topologie d'orbifold complexe quotient de celle de $\Teich(S,P)$. C'est aussi le quotient de l'espace de Torelli par le groupe $\MCG(S,P)/T(S,P) = \Diff^+(S,P) / \Diff_*^+(S,P)$.

La forme d'intersection algébrique est une forme symplectique non dégénérée sur $H^1(S,\Z)$, ce qui permet de définir un morphisme du groupe modulaire $\MCG(S,P)$ à valeurs dans le groupe symplectique $\Sp_{2g}(\Z)$. Le noyau de ce morphisme est le sous-groupe de Torelli $T(S,P)$, et ce morphisme est surjectif donc le groupe quotient $\MCG(S,P)/T(S,P)$ est isomorphe à $\Sp_{2g}(\Z)$ (voir~\cite{mcg}).

\bigskip

Soit $X$ une surface hyperbolique connexe, orientée, d'aire finie, de genre $g$, avec $q$ cuspides, marquée par un difféomorphisme préservant l'orientation $h : S \bs P \ra X$. Le premier groupe de cohomologie $H^1(X,\R)$ est un espace vectoriel réel, considérons le sous-espace vectoriel $H^1_c(X,\R)$ de $H^1(X,\R)$ constitué des classes de cohomologie de $1$-formes différentielles fermées à support compact. Considérons l'inclusion $\iota : S \bs P \ra S$, elle induit un isomorphisme $\iota^* : H^1(S,\R) \ra H^1_c(S \bs P,\R)$. Ainsi ${h^{-1}}^* \circ \iota^*$ est un isomorphisme entre $H^1(S,\R)$ et $H^1_c(X,\R)$, donc ce dernier est de dimension réelle $2g$. 

Cet espace vectoriel est également isomorphe au premier groupe de cohomologie $L^2$ réduite de $X$, ainsi qu'à l'espace vectoriel des $1$-formes différentielles harmoniques sur $X$ (dans le cas compact c'est le théorème de Hodge, dans le cas général voir~\cite{hodge_hyperbolic} et \cite[Corollary~1.6, p.~7 et Theorem~2.16, p.~27]{carron}). Nous utiliserons ici ces trois points de vue, et lorsqu'il faudra choisir un représentant d'une classe de cohomologie à support compact, nous choisirons l'unique représentant harmonique.

Sur l'espace $H^1_c(X,\R)$ vu comme premier groupe de cohomologie $L^2$ réduite, nous disposons du produit scalaire $L^2$, noté $\<\cdot,\cdot\>_X$. On peut également voir ce produit scalaire gr\^{a}ce à l'étoile de Hodge $*_X$. L'étoile de Hodge est une anti-involution de l'espace vectoriel des $1$-formes différentielles réelles sur $X$, définie comme la précomposition, sur chaque plan tangent, par la rotation d'un quart de tour dans le sens positif (c'est-à-dire la multiplication par $i$, pour la structure complexe sur $X$ correspondante). L'étoile de Hodge définit par restriction une involution sur l'espace des $1$-formes différentielles harmoniques $H^1_c(X,\R)$. Le produit scalaire sur $H^1_c(X,\R)$ s'exprime alors ainsi :
$$ \forall \omega, \omega' \in H^1_c(X,\R),\; \<\omega,\omega'\>_X = \int_X \omega \wedge *_X \omega'.$$
On notera $\|\cdot\|_X$ la norme euclidienne ($L^2$) associée sur $H^1_c(X,\R)$.

Par ailleurs, la forme d'intersection algébrique $\Int$ sur $H^1_c(X,\R)$ (image de la forme d'intersection algébrique sur $H^1(S,\R)$) est alternée et non dégénérée, et est donnée par 
$$ \forall \omega, \omega' \in H^1_c(X,\R),\; \Int(\omega,\omega') = \int_X \omega \wedge \omega' = \<\omega, -*_X \omega'\>_X.$$

\bigskip

Remarquons que l'espace vectoriel $H^1_c(X,\R)$, muni de ce produit scalaire et de la forme d'intersection, est isomorphe isométriquement symplectiquement à l'espace euclidien $\R^{2g}$, muni de la forme symplectique standard
$$ b(x,x') = \sum_{j=1}^g x_{j}x'_{g+j} - x_{g+j}x'_{j}.$$
Choisissons un tel isomorphisme linéaire isométrique symplectique
$$\varphi_X : H^1_c(X,\R) \ra \R^{2g}.$$
Le sous-groupe de $\GL_{2g}(\R)$ préservant la forme $b$ est $\U(b) = \Sp_{2g}(\R)$.

Considérons $\mathcal{E}^{b}_{2g}$ l'espace des réseaux de covolume $1$ de l'espace euclidien $\R^{2g}$, autoduaux pour la forme symplectique $b$, marqués par $\Z^{2g}$, à isométrie près. Il s'identifie à l'espace symétrique hermitien $\mathcal{E}^{b}_{2g} = \Sp_{2g}(\R) / \SU(g)$, d'après le lemme~\ref{lem:type_symplectique}. Fixons une base symplectique de $H^1(S,\Z)$, ce qui nous permettra d'identifier $H^1(S,\Z)$ avec $\Z^{2g}$.

\bigskip
%
L'application période de l'espace de Torelli à valeurs dans l'espace symétrique $\Sp_{2g}(\R) / \SU(g)$ ou son quotient de l'espace des modules dans l'espace localement symétrique $\Sp_{2g}(\Z) \bs \Sp_{2g}(\R) / \SU(g)$ est un objet très classique (voir par exemple~\cite{griffiths}, \cite{debarre}, \cite{mess}, \cite{mcg}). La nouvelle formulation ci-dessous facilitera notre compactification de l'espace de Torelli, et nous ferons le lien avec la définition classique très prochainement. Cette formulation a de plus l'avantage d'inclure le cas non compact (i.e. $P \neq \emptyset$).

Si $[X,h] \in \Tor(S,P)$, alors $\varphi_X \circ {h^{-1}}^* \circ \iota^* : H^1(S,\Z) \ra \R^{2g}$ est un morphisme de groupes. Définissons l'application période
\beq p : \Tor(S,P) & \ra & \mathcal{E}^b_{2g} \\
\ [X,h] & \mapsto & \left[ \varphi_X \circ {h^{-1}}^* \circ \iota^*|_{H^1(S,\Z)} \right]. \eeq

\blem \label{lem:periode_continue} L'application période est bien définie et continue. \elem

\bp
Puisque les trois applications linéaires $\iota^* : H^1(S,\R) \ra H^1_c(S \bs P,\R)$, ${h^{-1}}^* : H^1_c(S \bs P,\R) \ra H^1_c(X,\R)$ et $\varphi_X : H^1_c(X,\R) \ra \R^{2g}$ sont des isomorphismes linéaires, l'image $\varphi_X \circ {h^{-1}}^* \circ \iota^* (H^1(S,\Z))$ est un réseau de $\R^{2g}$.

Vérifions que ce réseau $\varphi_X \circ {h^{-1}}^* \circ \iota^*(H_1(S,\Z))$ de $\R^{2g}$ est bien autodual par rapport à $b$. Soit $c \in H^1(S,\R)$ tel que, pour tout $d \in H^1(S,\Z)$, nous ayons $b(\varphi_X \circ {h^{-1}}^* \circ \iota^*(c),\varphi_X \circ {h^{-1}}^* \circ \iota^*(d)) \in \Z$. Fixons $d \in H^1(S,\Z)$. Puisque $\varphi_X : H^1_c(X,\R) \ra \R^{2g}$ préserve les formes symplectiques, on sait que $\Int({h^{-1}}^* \circ \iota^*(c),{h^{-1}}^* \circ \iota^*(d)) \in \Z$. L'isomorphisme ${h^{-1}}^* : H^1_c(S \bs P,\R) \ra H^1_c(X,\R)$ provient d'un difféomorphisme préservant l'orientation $h^{-1} : X \ra S \bs P$, donc préserve les formes d'intersection, ainsi $\Int(\iota^*(c),\iota^*(d)) \in \Z$. Enfin, pour les formes différentielles $\iota^*(c)$ et $\iota^*(d)$ à support compact, leur nombre d'intersection sur $S \bs P$ et sur $S$ est le même, donc $\Int(c,d) \in \Z$. Or on sait que le réseau standard $\Z^{2g} = H^1(S,\Z)$ de $\R^{2g}$ est autodual pour la forme symplectique, donc $c \in H^1(S,\Z)$. Ainsi le réseau $\varphi_X \circ {h^{-1}}^* \circ \iota^*(H_1(S,\Z))$ de $\R^{2g}$ est autodual par rapport à $b$.

De plus, d'après le lemme~\ref{lem:type_symplectique}, il est dans l'orbite de $\Z^{2g}$ sous $\Sp_{2g}(\R)$, donc il est de covolume $1$.

Supposons que $[X',h']$ soit égal à $[X,h]$ dans l'espace de Torelli, alors il existe une isométrie préservant l'orientation $s : X \ra X'$ telle que ${h^{-1}}^* = s^* \circ {h'^{-1}}^* : H^1(S \bs P,\R) \ra H^1_c(X,\R)$. Or l'application $s^* : H^1_c(X',\R) \ra H^1_c(X,\R)$ est une isométrie linéaire préservant la forme d'intersection, donc les deux réseaux marqués $\varphi_X \circ {h^{-1}}^* \circ \iota^*|_{H^1(S,\Z)}$ et $\varphi_{X'} \circ {h'^{-1}}^* \circ \iota^*|_{H^1(S,\Z)}$ sont isométriques symplectiquement, et égaux dans $\mathcal{E}^b_{2g}$.

Remarquons également que l'application $p$ ne dépend pas du choix de l'isomorphisme isométrique symplectique $\varphi_X : H^1_c(X,\R) \ra \R^{2g}$ : deux tels isomorphismes diffèrent par composition au but d'une isométrie symplectique de $\R^{2g}$, ce qui donne une même image dans l'espace $\mathcal{E}^b_{2g}$.

\bigskip

Montrons que l'application $p$ est continue. Soit $([X_n,h_n])_{n \in \N}$ une suite convergeant vers $[X,h]$ dans $\Tor(S,P)$. Dans le fibré localement trivial donné par $H^1_c(\cdot,\R)$ au-dessus de l'espace de Torelli $\Tor(S,P)$, la suite $(H^1_c(X_n,\R))_{n}$ converge vers $H^1_c(X,\R)$. Ceci nous permet donc d'identifier chaque $H^1_c(X_n,\R)$ avec $H^1_c(X,\R)$.

La forme d'intersection sur ${h_n^{-1}}^* \circ \iota^*(H^1(S,\Z)) \subset H^1_c(X_n,\R)$ est donnée par l'image de la forme d'intersection sur $H^1(S,\Z)$. Puisque la suite $(X_n)_{n \in \N}$ converge vers $X$ dans l'espace des modules $\Mod(S,P)$, la suite des involutions de Hodge $(*_{X_n})_{n \in \N}$ converge vers l'involution de Hodge $*_X$ sur $H^1_c(X,\R)$. Ainsi la suite des produits scalaires $(\<\cdot,\cdot\>_{X_n})_{n \in \N}$ converge vers $\<\cdot,\cdot\>_X$.

Il reste à montre que la suite des marquages converge : on sait que la suite $(h_n)_{n \in \N}$ converge vers $h$ modulo $\Diff^+_{H^1}(S,P)$, donc la suite $({h_n^{-1}}^*)_{n \in \N}$ converge vers ${h^{-1}}^*$ pour la topologie des applications linéaires de $H^1_c(S \bs P,\R)$ dans $H^1_c(X_n,\R)$ identifié à $H^1_c(X,\R)$.

Nous avons montré que la suite $(p([X_n,h_n]))_{n \in \N}$ convergeait vers $p([X,h])$, l'application période est donc continue. \ep

\bthm[Griffiths] Dans le cas compact ($P=\emptyset$), l'application $p$ est holomorphe. L'adhérence de son image $\ov{p(\Tor(S))}$ est une sous-variété analytique complexe fermée de $\mathcal{E}^b_{2g}$, et l'image $p(\Tor(S))$ est le complémentaire dans son adhérence d'une sous-variété analytique complexe fermée. \ethm

\bp Rappelons la définition classique de l'application période (voir~\cite[Chapter~3, §2.1]{mcg}). Nommons $(\alpha_1,\ldots,\alpha_g,\beta_1,\ldots,\beta_g)$ la base de $H^1(S,\Z)$ qui s'envoie sur la base canonique de $\Z^{2g}$ par l'isomorphisme fixé entre $H^1(S,\Z)$ et $\Z^{2g}$. Considérons la base symplectique $(a_1,\ldots,a_g,b_1,\ldots,b_g)$ de $H_1(S,\Z)$ obtenue par dualité de Poincaré. Fixons $[X,h] \in \Tor(S)$. Soit $(\omega_1, \ldots, \omega_g, \omega'_1, \ldots, \omega'_g)$ une base de $H^1_c(X,\R)$, image de la base canonique de $\R^{2g}$ par un isomorphisme isométrique respectant la forme symplectique, et choisissons pour $\varphi_X$ cet isomorphisme.

Comme nous sommes ici dans le cas compact ($P = \emptyset$), l'application $h_*$ est définie de $H_1(S,\R)$ dans $H_1(X,\R)$. La matrice période de $[X,h]$ est alors définie comme la matrice par blocs
$$ \Omega = \left( \begin{array}{cc} \left(\int_{h_*(a_j)} \omega_k \right)_{j,k \in \lb 1,g \rb} & \left(\int_{h_*(a_j)} \omega'_k \right)_{j,k \in \lb 1,g \rb}  \\
\left(\int_{h_*(b_j)} \omega_k \right)_{j,k \in \lb 1,g \rb} &  \left(\int_{h_*(b_j)} \omega'_k \right)_{j,k \in \lb 1,g \rb} \end{array} \right) \in \M_{2g}(\R).$$
Puisque $(\omega_1, \ldots, \omega_g, \omega'_1, \ldots, \omega'_g)$ est une base symplectique de $H^1_c(X,\R)$, on en déduit que $\Omega \in \Sp_{2g}(\R)$. Cette matrice ne dépend que du choix de la base $(\omega_1, \ldots, \omega_g, \omega'_1, \ldots, \omega'_g)$, et deux choix de telles bases diffèrent par l'action d'une isométrie de $H^1_c(X,\R)$ préservant la forme d'intersection, c'est-à-dire par l'action de $\Sp_{2g}(\R) \cap \SO_{2g}(\R) = \SU(g)$. Ainsi la matrice période appartient à l'espace symétrique $\mathcal{E}^b_{2g} = \Sp_{2g}(\R) / \SU(g)$ : cela définit l'application période classique $\Tor(S) \ra \mathcal{E}^b_{2g}$.

\bigskip

Explicitons l'identification entre l'espace $\mathcal{E}^b_{2g}$ des réseaux marqués autoduaux de $\R^{2g}$ de covolume $1$ et l'espace symétrique $\Sp_{2g}(\R) / \SU(g)$, décrite au début de la partie~\ref{sec:thurston}, dans notre cas particulier. La classe d'isométrie d'un réseau marqué $f: \Z^{2g} \ra \R^{2g}$ autodual pour $b$ et de covolume $1$, s'identifie à la transposée de la matrice de l'isomorphisme linéaire $f$ de $\R^{2g}$, modulo multiplication par $\SU(g)$ à droite. Dans le cas de la classe du réseau marqué $p([X,h])$, la matrice ainsi obtenue est
$$ \Omega' = \left( \begin{array}{cc} \left(\<{h^{-1}}^*(\alpha_j),\omega_k\> \right)_{j,k \in \lb 1,g \rb} & \left(\<{h^{-1}}^*(\alpha_j),\omega'_k\> \right)_{j,k \in \lb 1,g \rb}  \\
\left(\<{h^{-1}}^*(\beta_j),\omega_k\> \right)_{j,k \in \lb 1,g \rb} &  \left(\<{h^{-1}}^*(\beta_j),\omega'_k\> \right)_{j,k \in \lb 1,g \rb} \end{array} \right) \in \M_{2g}(\R).$$
Or chaque produit scalaire est égal à une intégrale intervenant dans la matrice $\Omega$ : par exemple, pour tous $j,k \in \lb 1,g \rb$, nous avons
$$ \<{h^{-1}}^*(\alpha_j),\omega_k\> = \Int({h^{-1}}^*(\alpha_j), *_X \omega_k) = \int_{h_*(a_j)} *_X \omega_k = \int_{h_*(a_j)} \omega'_k,$$
par dualité de Poincaré. Dans les autres cas il faut utiliser la relation $*_X \omega'_k = -\omega_k$, et on remarque alors que $\Omega'=\Omega J'$, où $J'$ est la matrice par blocs
$$ J' = \left( \begin{array}{cc} 0 & -I_g \\ I_g & 0 \end{array} \right) \in \SU(g).$$
Ainsi les deux définitions de l'application période sont les mêmes.

\bigskip

D'après~\cite[Proposition~9.3, p.~156]{griffiths}, l'application période est holomorphe. De plus, d'après~\cite[Theorem~9.6]{griffiths}, l'adhérence de l'image de l'application période est une sous-variété analytique complexe fermée, qui contient $\Im p$ comme complémentaire d'une sous-variété analytique complexe fermée. \ep

\bthm[Torelli] \label{thm:torelli} Dans le cas compact ($P=\emptyset$), l'application $p$ est un revêtement ramifié de degré deux sur son image, ramifié sur le lieu hyperelliptique de $\Tor(S)$. \ethm

\bp La forme forte du théorème de Torelli, formulée dans l'article~\cite[Theorem~1, p.~783]{mess}, énonce que l'application période est un revêtement ramifié de degré deux sur son image, ramifié sur le lieu hyperelliptique de $\Tor(S)$. \ep

Considérons l'application (introduite en introduction dans le cas $P =\emptyset$)
\beq \psi : \Tor(S,P) & \ra & {\R_+}^{\!H^1(S,\Z)} \\
\ [X,h] & \mapsto & \left(c \mapsto  \| {h^{-1}}^* \circ \iota^*(c) \|_X \right), \eeq
et son image après composition par l'application canonique ${\R_+}^{\!H^1(S,\Z)} \ra \P({\R_+}^{\!H^1(S,\Z)})$, c'est-à-dire
\beq \ov{\psi} : \Tor(S,P) & \ra & \P({\R_+}^{\!H^1(S,\Z)}) \\
\ [X,h] & \mapsto & \left[c \mapsto  \| {h^{-1}}^* \circ \iota^*(c) \|_X \right]. \eeq
Les applications $\psi$ et $\ov{\psi}$ sont $\Sp_{2g}(\Z)$-équivariantes pour les actions de $\Sp_{2g}(\Z)$ (s'identifiant à $\MCG(S,P) / T(S,P)$) sur $\Tor(S,P)$ et de $\Sp_{2g}(\Z)$ sur $\Z^{2g} = H^1(S,\Z)$. 

L'application $\ov{\psi}$ est égale à la composition $\phi \circ p$, où $\phi :\mathcal{E}^b_{2g} \ra \P({\R_+}^{\!H^1(S,\Z)})$ est le plongement de la proposition~\ref{pro:plongement}.

Dans le cas compact ($P = \emptyset$), d'après le théorème~\ref{thm:torelli}, l'application $\ov{\psi}$ est un revêtement ramifié d'ordre deux sur son image. Dans le cas non compact, d'après le lemme~\ref{lem:periode_continue}, l'application $\ov{\psi}$ est continue.

Notons $\Tor(S,P) \cup \{\infty\}$ la compactification d'Alexandrov de $\Tor(S,P)$, et considérons le plongement diagonal de $\Tor(S,P)$ dans
$$\P(\R_+^{H_1(S,\Z)}) \times (\Tor(S,P) \cup \{\infty\})$$
donné par le produit de $\ov{\psi}$ et de l'inclusion. Appelons \df{compactification de Thurston} de l'espace de Torelli, notée $\ov{\Tor(S,P)}^\mathcal{T}$, l'adhérence de l'image de ce plongement. C'est une $\Sp_{2g}(\Z)$-compactification.

D'autre part, considérons le plongement diagonal de $\Tor(S)$ dans $\ov{\mathcal{E}^b_{2g}}^\mathcal{S} \times (\Tor(S) \cup \{\infty\})$ donné par le produit de $p$ et de l'inclusion. Appelons \df{compactification de Satake} de l'espace de Torelli, notée $\ov{\Tor(S)}^\mathcal{S}$, l'adhérence de l'image de ce plongement. C'est une $\Sp_{2g}(\Z)$-compactification.

\bthm \label{thm:torelli_isomorphes} Les compactifications de Thurston $\ov{\Tor(S,P)}^\mathcal{T}$ et de Satake $\ov{\Tor(S,P)}^\mathcal{S}$ de l'espace de Torelli sont $\Sp_{2g}(\Z)$-isomorphes. \ethm

\bp Il suffit de remarquer que les deux compactifications de $p(\Tor(S,P)) \subset \mathcal{E}^b_{2g}$, où $p(\Tor(S,P))$ est plongé dans les deux compactifications $\ov{\mathcal{E}^b_{2g}}^\mathcal{T}$ et $\ov{\mathcal{E}^b_{2g}}^\mathcal{S}$, sont $\Sp_{2g}(\Z)$-isomorphes : c'est une conséquence du théorème~\ref{thm:thurston_satake}. \ep

\section{Stratification d'une partie du bord de l'espace de Torelli}

Une stratification naturelle du bord d'une compactification permet de mieux en comprendre la structure, comme dans les travaux de Harvey sur l'espace des modules (voir~\cite{harvey_geometric} et \cite{harvey_boundary}). Nous allons décrire une stratification naturelle d'une partie du bord de la compactification de Thurston de l'espace de Torelli.

Dans toute cette partie, nous fixerons comme au début la partie~\ref{sec:torelli} une surface $S$ lisse compacte connexe orientée de genre $g$, et $P$ une partie finie fixée de $S$ de cardinal $q$, telle que $2g-2+q > 0$.

\subsection{Comparaison entre norme stable et norme euclidienne}

Fixons une surface hyperbolique $X$ d'aire finie connexe orientée. Appelons \df{systole séparante} (resp. non séparante) de $X$ la borne inférieure (qui est atteinte) des longueurs des courbes fermées simples de $X$ non homotopes à zéro, qui ne bordent pas une unique cuspide, et séparantes (resp. non séparantes). Notons $\sns(X)$ la systole non séparante de $X$.

Pour les rappels qui suivent, nous renvoyons à la thèse de Daniel Massart (\cite{massart_these}).

Définissons une norme naturelle sur $H_1(X,\R)$, la \df{norme stable} (voir \cite{gromov} et \cite{massart}). Si $\gamma$ est une courbe fermée simple non homotope à zéro sur $X$, notons $\lg(\gamma)$ la longueur hyperbolique de l'unique géodésique dans la classe d'homotopie libre de $\gamma$. Soit $c \in H_1(X,\R)$. Alors on peut représenter $c$ par un $1$-cycle singulier de la forme $\sum_{i \in I} x_i \alpha_i$, où $I$ est fini, $x_i \in \R$ (et $x_i \in \Z$ si $c \in H_1(X,\Z)$) et $\alpha_i$ est une courbe fermée simple non homologue à zéro sur $X$, qu'on appelle \df{multicourbe}. On appelle \df{support} d'une telle multicourbe la réunion des courbes $\alpha_i$ pour lesquelles $x_i \neq 0$. Définissons alors la \df{norme stable} de $c$ par
$$  \| c \|_s = \inf_{c = [\sum_{i \in I} x_i \gamma_i]} \sum |x_i| \lg(\gamma_i),$$
où $\sum_{i \in I}^n x_i \gamma_i$ parcourt toutes les multicourbes dans la classe d'homologie de $c$.

Notons $\theta_X : H_1(X,\R) \ra H^1(X,\R)$ l'isomorphisme de dualité de Poincaré. Définissons la norme stable sur sur $H^1(X,\R)$ par $\|\omega\|_{X,1} = \| \theta_X^{-1}(\omega)\|_s$, où $\omega \in H^1(X,\R)$. Nous avons par ailleurs défini une norme euclidienne sur $H^1_c(X,\R)$, notons-la $\|\cdot\|_{X,2}$.

Dans sa thèse (\cite{massart_these}), Daniel Massart compare les deux normes précédemment introduites.

\bthm[Massart] \label{thm:l1l2} Il existe deux constantes strictement positives $a$ et $b$, ne dépendant que du genre de $X$, telles que pour tout $\omega \in H^1_c(X,\R)$ nous ayons
$$ \hspace{2cm} a\; \| \omega \|_{X,1} \leq  \| \omega \|_{X,2} \leq \frac{b}{\sns(X)^2}\;  \| \omega \|_{X,1}.$$
\ethm

\bp Dans la partie~4.2 de~\cite{massart_these}, ce théorème est énoncé pour une surface hyperbolique compacte sans bord. Pour la première inégalité, D.~Massart utilise la densité des formes de Strebel dans l'ensemble des $1$-formes différentielles harmoniques, et ce résultat est encore vrai pour une surface de volume fini (voir~\cite{douady_hubbard}). Et la deuxième inégalité n'utilise même pas l'hypothèse de courbure, et elle reste vraie pour une surface de volume fini.
\ep

Nous allons nous servir de ce théorème pour étudier l'image de $\Tor(S,P)$ par l'application $\psi$.

\blem \label{lem:sns_minoree} Soit $([X_n,h_n])_{n \in \N}$ une suite dans $\Tor(S,P)$ telle que la suite $(\psi([X_n,h_n]))_{n \in \N}$ soit bornée dans ${\R_+}^{\!H^1(S,\Z)}$. Alors la systole non séparante $(\sns(X_n))_{n \in \N}$ est minorée par une constante strictement positive. \elem

\bp
Par contraposée, supposons que la systole non séparante de $(X_n)_{n \in \N}$ tende vers $0$ : pour tout $n \in \N$, soit donc $\gamma_n$ une géodésique fermée simple de $X_n$ non séparante et non cuspidale, dont la longueur $\lg(\gamma_n)=\sns(X_n)$ tend vers zéro. La courbe $\gamma_n$ est non triviale et non cuspidale, donc a un nombre d'intersection non nul avec l'un des générateurs de
$$\theta_{X_n}^{-1}(H^1_c(X_n,\Z)) = {h_n}_* \circ \theta_{S \bs P}^{-1} \circ \iota^*(H^1(S,\Z)),$$
images des éléments de la base choisie de $H^1(S,\Z)$. Quitte à extraire, on peut supposer qu'il existe $\omega \in H^1(S,\Z)$ tel que le nombre d'intersection algébrique de $\gamma_n$ avec ${h_n}_* \circ \theta_{S \bs P}^{-1} \circ \iota^*(\omega)$ soit non nul pour tout $n$. D'après~\cite[Corollary~4.1.2, p.~95]{buser}, on en déduit que la longueur d'un représentant minimal de ${h_n}_* \circ \theta_{S \bs P}^{-1} \circ \iota^*(\omega)$ est minorée par $2\argsh\left((\sh(\frac{\lg(\gamma_n)}{2}))^{-1}\right)$, ce qui tend vers l'infini lorsque $n$ tend vers $+\infty$. Ainsi la norme stable $\| {h_n}_* \circ \theta_{S \bs P}^{-1} \circ \iota^*(\omega) \|_s$ tend vers l'infini lorsque $n$ tend vers $+\infty$.

Or, d'après le théorème~\ref{thm:l1l2}, on sait que
\beq
\| {h_n^{-1}}^* \circ \iota^*(\omega) \|_{X_n,2} & \geq & a \| {h_n^{-1}}^* \circ \iota^*(\omega) \|_{X_n,1} \\
 & \geq & a \| {h_n^{-1}}^* \circ \iota^*(\omega) \|_{X_n,1} \\
 & \geq & a \| {h_n}_* \circ \theta_{S \bs P}^{-1} \circ \iota^*(\omega) \|_s \ral{n \ra +\infty} +\infty, \eeq
donc la suite de fonctions de longueur $(\psi([X_n,h_n]))_{n \in \N}$ n'est pas bornée.
\ep

\blem \label{lem:l_propre} Soit $([X_n,h_n])_{n \in \N}$ une suite dans $\Tor(S,P)$ telle que la suite $(\psi([X_n,h_n]))_{n \in \N}$ converge vers $\ell$ dans ${\R_+}^{\!H^1(S,\Z)}$. Alors la suite $(p([X_n,h_n]))_{n \in \N}$ converge dans $\mathcal{E}^b_{2g}$, et en particulier la fonction de longueur $\ell$ est propre. \elem

\bp
D'après le lemme~\ref{lem:sns_minoree}, on sait que la systole non séparante $\sns(X_n)$ est minorée par une constante $d>0$. Ainsi d'après le théorème~\ref{thm:l1l2}, pour tout $\omega \in H^1(S,\Z)$ non nul et tout $n \in \N$, nous avons
$$\psi([X_n,h_n])(\omega) =\|{h_n^{-1}}^* \circ \iota^*(\omega)\|_{X_n,2} \geq a \|{h_n^{-1}}^* \circ \iota^*(\omega) \geq ad >0.$$
Ainsi la systole du réseau $\varphi_{X_n} \circ {{h_n}^{-1}}^*(H^1(S,\Z))$ de $\R^{2g}$ de covolume $1$ est minorée par $ad$, donc d'après le critère de Mahler (voir~\cite{mahler}) on peut supposer quitte à extraire que cette suite de réseaux de $\R^{2g}$ non marqués converge vers un réseau $\Lambda \subset \R^{2g}$ de covolume $1$ et de systole minorée par $ad$. Ce réseau $\Lambda$ est autodual pour la forme symplectique $b$ standard de $\R^{2g}$. Puisque la suite de fonctions de longueur $(\psi([X_n,h_n]))_{n \in \N}$ converge, on en déduit quitte à extraire que la suite de réseaux marqués $(\varphi_{X_n} \circ {{h_n}^{-1}}^*|_{H^1(S,\Z)})_{n \in \N}$ converge vers le réseau marqué $(\Lambda,h)$, où $h : H^1(S,\Z) \ra \Lambda$ est un marquage. Ainsi la suite $(p([X_n,h_n]))_{n \in \N}$ converge vers $[\Lambda,h]$ dans $\mathcal{E}^b_{2g}$.

En particulier, la suite de fonctions de longueur $(\psi([X_n,h_n]))_{n \in \N}$ converge vers $\ell = \phi([\Lambda,h])$, qui est une application propre.
\ep

\subsection{Le complexe des courbes séparantes}

Pour plus de détails sur cette partie, on pourra consulter~\cite[Chapter~3, §~2.4]{mcg}. Comme au début la partie~\ref{sec:torelli}, on considère une surface $S$ lisse compacte connexe orientée de genre $g$, et $P$ une partie finie fixée de $S$ de cardinal $q$. On dit qu'une courbe fermée simple $\gamma$ dans $S \bs P$ est \df{séparante} si $S \bs (P \cup \gamma)$ n'est pas connexe. On dit que $\gamma$ est \df{cuspidale} si elle bord un disque épointé dans $S \bs P$, et \df{triviale} si elle borde un disque dans $S \bs P$.

Définissons le \df{complexe des courbes séparantes} $K^{sep}(S,P)$ de $S$. Les sommets de $K^{sep}(S,P)$ sont les classes d'isotopie de courbes fermées simples séparantes de $S \bs P$ non cuspidales et non triviales. Les classes d'isotopie de telles courbes $\gamma_1, \ldots ,\gamma_k$ forment un ($k$-$1$)-simplexe de $K^{sep}(S,P)$ si ces courbes sont disjointes et si leurs classes d'isotopie sont distinctes.

Notons $\Sigma K^{sep}(S,P)$ l'ensemble des simplexes de $K^{sep}(S,P)$.

\subsection{Description des strates}

\label{subsec:torelli_strates}

\begin{figure}[!h]
\begin{center}
\includegraphics[height=4cm]{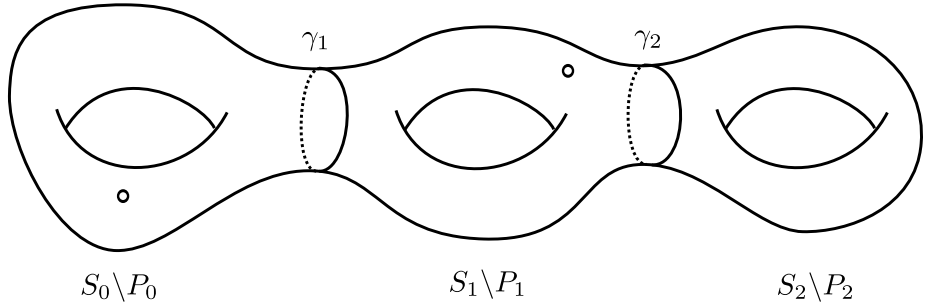}
\caption{Découpage de la surface $S$}
\label{fig:decoupe}
\end{center}
\end{figure}

Soit $\sigma=\{[\gamma_1],\ldots,[\gamma_k]\}$ un ($k$-$1$)-simplexe de $K^{sep}(S,P)$. Alors $S \bs (P \cup \bigcup_{j=1}^k \gamma_j)$ a $k+1$ composantes connexes. \`{A} difféomorphisme préservant l'orientation près (fixé), ces composantes s'écrivent $(S_j \bs P_j)_{j \in \lb 0,k\rb}$, où $S_j$ est une surface lisse compacte connexe orientée de genre $g_j$ et $P_j \subset S_j$ est un ensemble de cardinal $q_j$ fini, tels que $2-2g_j-q-j<0$. De plus, on a les égalités $\sum_{j=0}^{k} g_j = g$ et $\sum_{j=0}^{k} q_j = 2k+q$. Pour tout $j \in \lb 0,k \rb$, notons l'inclusion $\iota_j : S_j \bs P_j \ra S_j$. Voir la figure~\ref{fig:decoupe}.

\bigskip

Considérons, pour tout $j \in \lb 0,k\rb$, l'application $\kappa_j : S \ra S_j$ qui écrase chaque $S_{j'}$, pour $j' \neq j$. Considérons alors l'application
\beq \xi_\sigma : \bigoplus_{j=0}^k H^1(S_j,\R) & \ra & H^1(S,\R) \\
(\omega_j)_{j \in \lb 0,k \rb} & \ra & \sum_{j=0}^k {\kappa_j}^*(\omega_j),\eeq
elle s'appelle le \df{scindement symplectique} de $H^1(S,\R)$ associé à $\sigma$. Remarquons qu'elle définit aussi un scindement symplectique de $\bigoplus_{j=0}^k H^1(S_j,\Z) \ra H^1(S,\Z)$, encore noté $\xi_\sigma$. Il est dit symplectique car cette décomposition est orthogonale pour la forme d'intersection.

\bigskip

Considérons, pour tout $j \in \lb 0,k \rb$, l'application $\psi_j$ de $\Tor(S_j,P_j)$ à valeurs dans ${\R_+}^{\!H^1(S_j,\Z)}$, définie pour le théorème~\ref{thm:torelli_isomorphes}.

\bigskip

Appelons \df{strate} au bord de l'espace de Torelli de $S$ associée au simplexe $\sigma$ la variété différentielle produit $\Tor_\sigma(S,P) = \prod_{j=0}^k \Tor(S_j,P_j)$. Nous allons la relier à l'espace de Torelli total $\Tor(S,P)$ grâce à l'application
\beq \psi_\sigma : \Tor_\sigma(S,P) & \ra &  {\R_+}^{\!H^1(S,\Z)} \\
\left([X_j,h_j]\right)_{j \in \lb 0,k \rb} & \mapsto & \left\{ \omega \mapsto \left( \sum_{j=0}^k \psi_j([X_j,h_j]) \circ ({\xi_\sigma}^{-1})_j(\omega)^2 \right)^\frac{1}{2} \right\}.\eeq

\bthm L'adhérence de l'image de l'application $\psi$ dans l'espace ${\R_+}^{\!H^1(S,\Z)}$ est la réunion disjointe 
$$ \ov{\psi(\Tor(S,P))} = \psi(\Tor(S,P)) \sqcup \bigsqcup_{\sigma \in \Sigma K^{sep}(S,P)} \psi_\sigma(\Tor_\sigma(S,P)).$$ \ethm

\bp
Soit $([X_n,h_n])_{n \in \N}$ une suite de $\Tor(S,P)$ telle que la suite $(\psi([X_n,h_n]))_{n \in \N}$ converge vers $\ell \in \partial \psi(\Tor(S,P))$. D'après le lemme~\ref{lem:sns_minoree}, on sait que la systole non séparante $(\sns(X_n))_{n \in \N}$ est minorée par une constante strictement positive.

Supposons que la systole de la suite $(X_n)_{n \in \N}$ soit minorée par une con\-stante strictement positive, alors d'après~\cite[Corollary~I.3.1.7, p.~63]{ceg} par exemple, on sait que la suite $(X_n)_{n \in \N}$ est bornée dans l'espace des modules $\Mod(S)$. Quitte à extraire, on peut donc supposer que la suite $(X_n)_{n \in \N}$ converge vers une surface hyperbolique $X$ dans $\Mod(S)$. Puisque la suite $([X_n,h_n])_{n \in \N}$ part à l'infini dans $\Tor(S,P)$, on en déduit que c'est la suite de marquages $(h_n)_{n \in \N}$ qui part à l'infini : il existe donc $\omega \in H_1(S,\Z)$ tel que la suite $({h_n^{-1}}^* \circ \iota^*(\omega))_{n \in \N}$ parte à l'infini, donc telle que la norme stable $(\|{h_n^{-1}}^* \circ \iota^*(\omega)\|_{X_n,1})_{n \in \N}$ tende vers l'infini (en effet, la surface hyperbolique $X$ a un nombre fini de géodésiques fermées simples de longueur majorée). En conséquence, d'après le théorème~\ref{thm:l1l2}, on en déduit que la norme euclidienne $(\|{h_n^{-1}}^* \circ \iota^*(\omega)\|_{X_n,2})_{n \in \N}$ tend vers l'infini, et donc la suite de fonctions de longueur $(\psi([X_n,h_n]))_{n \in \N}$ ne peut converger.

\bigskip

On en déduit donc que la systole séparante de la suite $(X_n)_{n \in \N}$ tend vers zéro. Fixons $\delta>0$ suffisamment petit, inférieur à la constante de Margulis du plan hyperbolique (voir par exemple~\cite[Theorem~I.2.2.2, p.~50]{ceg}), et inférieur aux systoles non séparantes des surfaces $(X_n)_{n \in \N}$. Considérons ${h_n}_*(\sigma_n)$ l'ensemble des classes d'homologie de géodésiques fermées simples séparantes de $X_n$ de longueur au plus $\delta$. Ces géodésiques sont disjointes, donc $\sigma_n$ est un simplexe du complexe $K^{sep}(S,P)$ des courbes séparantes. Quitte à extraire et quitte à diminuer $\delta$, on peut supposer que $\sigma_n$ est un $(k-1)$-simplexe dont la longueur de chacune des courbes tend vers zéro.

Quitte à précomposer les marquages $h_n : S \bs P \ra X_n$ par des éléments de $\Diff^+_{H^1}(S,P)$, comme il n'y a qu'un nombre fini d'orbites de simplexes de $K^{sep}(S,P)$ sous l'action de $\Diff^+_{H^1}(S,P)$, on peut supposer que le simplexe $(\sigma_n)_{n \in \N}$ est constant, égal à $\sigma$. Choisissons un représentant $\tilde{\sigma}$ de $\sigma$.

\bigskip

Notons $(S_0 \bs P_0), \ldots, (S_k \bs P_k)$ les composantes connexes de $S \bs \tilde{\sigma}$. Fixons un indice $j \in \lb 0,k \rb$, et notons $X_{n,j}$ la $j$ème composante connexe de $X_n \bs h_n(\sigma)$, où l'on a choisi les représentants géodésiques des courbes de $h_n(\sigma)$ : c'est une surface hyperbolique non complète d'aire finie. Quitte à changer le marquage $h_n : S \bs P \ra X$ par un difféoméorphisme isotope à l'identité, on peut supposer que $h_n$ envoie les courbes de $\tilde{\sigma}$ sur les représentants géodésiques de $h_n(\sigma)$. Ainsi, l'application $h_{n,j} = h_n|_{S_j \bs P_j} : S_j \bs P_j \ra X_{n,j}$ est un difféomorphisme. Alors, pour n'importe quel pointage de $X_{n,j}$ dans sa partie épaisse, la suite de surfaces hyperboliques ouvertes pointées $(X_{n,j})_{n \in \N}$ converge vers une surface hyperbolique $X_{\infty,j}$ pour la topologie de Gromov-Hausdorff, quitte à extraire. En effet, la systole de $X_{n,j}$ est minorée par $\delta$. Et puisque la longueur des courbes au bord de $X_{n,j}$ tend vers $0$, la surface $X_{\infty,j}$ est complète. De plus, quitte à changer les marquages $h_{n,j} : S_j \bs P_j \ra X_{n,j}$ par des difféomorphismes isotopes à l'identité, on peut supposer que la suite de marquages $(h_{n,j})_{n \in \N}$ converge vers un difféomorphisme $h_{\infty,j} : S_j \bs P_j \ra X_{\infty,j}$.

\bigskip

Fixons $j \in \lb 0,k \rb$, et $\omega_j \in H^1(S_j,\Z)$. Montrons que
$$ \| {h_n^{-1}}^* \circ \iota^* \circ \kappa_j^*(\omega_j) \|_{X_n,2} \ral{n \ra +\infty} \| {h_{\infty,j}^{-1}}^* \circ \iota^* \circ \kappa_j^*(\omega_j) \|_{X_{\infty,j},2}.$$
Pour tout $n \in \N$, notons $\eta_n = {h_n^{-1}}^* \circ \iota^* \circ \kappa_j^*(\omega_j) \in H^1_c(X_n,\R)$. Rappelons-nous que par convention $\eta_n$ désigne également l'unique $1$-forme différentielle harmonique dans cette classe de cohomologie. Puisque $\theta_{X_n}(\eta_n)$ peut être représenté par l'image par $h_n$ d'une multicourbe $c$ incluse dans $S_j \bs P_j$, et que la systole de la surface $X_{n,j} = h_n(S_j \bs P_j)$ est minorée par $\delta$, d'après le théorème~\ref{thm:l1l2}, on sait que la norme euclidienne de $\eta_n$ est majorée par $\|\eta_n\|_{X_n,2} \leq \frac{b}{\delta \log(\delta)} \|\eta_n\|_{X_n,1}$. Or
la multicourbe $c$ est la somme formelle d'un nombre fini de courbes fermées simples $c_i$ de $S_j \bs P_j$, or puisque la suite de marquages $(h_{n,j})_{n \in \N}$ converge vers $h_{\infty,j} : S_j \bs P_j \ra X_{\infty,j}$, la suite des longueurs $(\lg(h_{n,j}(c_i)))_{n \in \N}$ converge vers $\lg(h_{\infty,j}(c_i))$, donc est bornée. Ainsi la norme stable $\|\eta_n\|_{X_n,1}$ de $\eta_n$ est bornée, donc la norme euclidienne $\|\eta_n\|_{X_n,2}$ également.

\bigskip

Ainsi, pour tout $j' \in \lb 0,k \rb$, la restriction $\eta_{n,j'}$ de $\eta_n$ à $X_{n,j'}$ converge quitte à extraire vers une $1$-forme différentielle harmonique $\eta_{\infty,j'}$ sur $X_{\infty,j'}$. Si $j' \neq j$, alors $\eta_{n,j'}$ est cohomologue à zéro pour tout $n \in \N$, donc $\eta_{\infty,j'}$ est cohomologue à zéro. Puisque la surface est complète, on en déduit que $\eta_{\infty,j'}=0$. Et si $j'=j$, alors $\eta_{n,j}$ est cohomologue à ${h_{n,j}^{-1}}^* \circ \iota^* \circ \kappa_j^*(\omega_j)$, donc on en déduit que $\eta_{\infty,j}$ est cohomologue à ${h_{\infty,j}^{-1}}^* \circ \iota^* \circ \kappa_j^*(\omega_j)$. Par conséquent, pour les normes $L^2$, la suite $\|\eta_n\|_{X_n,2}^2 = \sum_{j'=0}^k \|\eta_{n,j'}\|^2$ converge, lorsque $n$ tend vers l'infini, vers
$$\sum_{j'=0}^k \|\eta_{\infty,j'}\|^2 = \|\eta_{\infty,j}\|^2 = \|{h_{\infty,j}^{-1}}^* \circ \iota^* \circ \kappa_j^*(\omega_j)\|_{X_{\infty,j},2}^2.$$

Nous avons donc montré que la suite $(\psi([X_n,h_n])(\kappa_j^*(\omega_j)))_{n \in \N}$ convergeait vers $\psi_j([X_{\infty,j},h_{\infty,j}])(\omega_j)$. Ceci prouve donc que la suite $(\psi([X_n,h_n]))_{n \in \N}$ converge dans ${\R_+}^{\!H^1(S,\Z)}$ vers
$$\psi_\sigma([X_{\infty,j},h_{\infty,j}]_{j \in \lb 0,k \rb}) \in \psi_\sigma(\Tor_\sigma(S,P)).$$

\bigskip

Réciproquement, montrons que $\psi_\sigma(\Tor_\sigma(S)) \subset \ov{\psi(\Tor(S))}$ pour tout $(k-1)$-simplexe $\sigma$ de $K^{sep}(S)$ : soit $([X_j,h_j])_{j \in \lb 0,k \rb} \in \Tor_\sigma(S)$. Fixons un entier $j \in \lb 0,k \rb$, les points marqués $P_j$ de $S_j$ sont de deux types : ceux qui correspondent à un point marqué de $P \subset S$ (nous les noterons $P_j^1$), et ceux qui correspondent à une courbe fermée simple de $\tilde{\sigma}$ (nous les noterons $P_j^2$).

Il existe une suite de surfaces hyperboliques $(X_{n,j})_{n \in \N}$ complètes d'aire finie avec $\Card P_j^2$ cuspides, et avec $\Card P_j^1$ composantes de bord totalement géodésique, chacune ayant pour longueur $\frac{1}{n+1}$, marquées par $h_{n,j} : S_j \bs P_j \ra X_{n,j} \bs \partial X_{n,j}$, telle que la suite de surfaces marquées $([X_{n,j},h_{n,j}])_{n \in \N}$ converge vers la surface marquée $[X_j,h_j]$ au sens suivant : pour tout compact $K$ de $S_j \bs P_j$, la suite d'espaces métriques compacts $(h_{n,j}(K))_{n \in \N}$ converge vers $h_j(K)$ au sens de Gromov-Hausdorff.

Pour tout $j \in \lb 1,k\rb$, la courbe $\gamma_j$ de $\tilde{\sigma}$ borde les surfaces $S_{j-1} \bs P_{j-1}$ et $S_j \bs P_j$. Pour tout $n \in \N$, il existe une composante de bord de $X_{n,j-1}$ et une de $X_{n,j}$ qui correspondent à $\gamma_j$ : puisque ces deux composantes ont la même longueur $\frac{1}{n+1}$, on peut recoller ces deux composantes ensemble. Notons $X_n$ le résultat du recollement des $k+1$ surfaces $X_{n,0},\ldots,X_{n,k}$, et $h_n : S \bs P \ra X_n$ un marquage qui coïncide avec chacun des $h_{n,j} : S_j \bs P_j \ra X_{n,j} \bs \partial X_{n,j}$ hors d'un petit voisinage de $\tilde{\sigma}$.

Alors d'après l'étude qui précède, la suite $(\psi([X_n,h_n]))_{n \in \N}$ converge vers $\psi_\sigma([X_j,h_j]_{j \in \lb 0,k \rb})$ dans ${\R_+}^{\!H^1(S,\Z)}$.

\bigskip

Enfin, montrons que les strates $\psi_\sigma(\Tor_\sigma(S,P))$ sont disjointes. Si $\ell \in {\R_+}^{\!H^1(S,\Z)}$ et $\sigma \in \Sigma K^{sep}(S,P)$, on dit que $\ell$ \df{se scinde selon} $\xi_\sigma$ si :
$$ \forall (\omega_0,\ldots,\omega_k) \in \prod_{j=0}^k H^1(S_j,\Z), \ell(\xi_\sigma(\omega_0,\ldots,\omega_k))^2 = \sum_{j=0}^k \ell(\xi_\sigma(0,\ldots,\omega_j,\ldots,0))^2.$$
On dit que $\ell$ se scinde s'il existe un simplexe $\sigma \in \Sigma K^{sep}(S,P)$ tel que $\ell$ se scinde selon $\xi_{\sigma}$.

Si $\ell \in \psi_\sigma(\Tor_\sigma(S))$, alors il est immédiat $\ell$ se scinde selon $\xi_\sigma$. Montrons que $\ell$ ne se scinde pas plus finement que $\xi_\sigma$, c'est-à-dire qu'il n'existe aucun simplexe $\sigma' \in \Sigma K^{sep}(S,P)$ contenant strictement $\sigma$ tel que $\ell$ se scinde selon $\xi_{\sigma'}$. Cela revient à montrer que, pour tout simplexe $\sigma \in \Sigma K^{sep}(S,P)$, aucune fonction de longueur $\ell \in \psi(\Tor(S,P))$ ne se scinde selon $\xi_\sigma$ : c'est le contenu du lemme~\ref{lem:non_scinde}.
\ep

\blem \label{lem:non_scinde} Aucune fonction de longueur de $\psi(\Tor(S,P))$ ne se scinde. \elem

\bp Supposons qu'il existe $[X,h] \in \Tor(S,P)$ tel que la fonction de longueur $\ell=\psi([X,h])$ se scinde selon $\xi_\sigma$, où $\sigma=\{[\gamma]\}$ et $\gamma$ est une courbe fermée simple séparante et non cuspidale de $S \bs P$ telle que $h(\gamma)$ soit géodésique.

%

\bigskip

Notons $S_1 \bs P_1$ et $S_2 \bs P_2$ les deux composantes connexes de $S \bs (P \cup \gamma)$.
%
%
Montrons que ${h^{-1}}^* \circ \xi_\sigma^{-1}(H^1(S_1,\R))$ et ${h^{-1}}^* \circ \xi_\sigma^{-1}(H^1(S_2,\R))$ sont orthogonaux pour le produit scalaire de $H^1_c(X,\R)$ : soient $\omega_1 \in H^1(S_1,\R)$ et $\omega_2 \in H^1(S_2,\R)$. Alors $\ell(\xi_\sigma^{-1}(\omega_1+\omega_2))^2 = \|{h^{-1}}^* \circ \xi_\sigma^{-1}(\omega_1+\omega_2)\|_X^2$ est égal à
$$ \|{h^{-1}}^* \circ \xi_\sigma^{-1}(\omega_1)\|_X^2 + 2 \<{h^{-1}}^* \circ \xi_\sigma^{-1}(\omega_1),{h^{-1}}^* \circ \xi_\sigma^{-1}(\omega_2)\>_X + \|{h^{-1}}^* \circ \xi_\sigma^{-1}(\omega_2)\|_X^2,$$
mais aussi à
$$\ell({h^{-1}}^* \circ \xi_\sigma^{-1}(\omega_1))^2 + \ell({h^{-1}}^* \circ \xi_\sigma^{-1}(\omega_1))^2 = \|{h^{-1}}^* \circ \xi_\sigma^{-1}(\omega_1)\|_X^2 + \|{h^{-1}}^* \circ \xi_\sigma^{-1}(\omega_2)\|_X^2,$$
ainsi $\<{h^{-1}}^* \circ \xi_\sigma^{-1}(\omega_1),{h^{-1}}^* \circ \xi_\sigma^{-1}(\omega_2)\>_X = 0$.

\bigskip

On en déduit que l'étoile de Hodge $*_X$ préserve ${h^{-1}}^* \circ \xi_\sigma^{-1}(H^1(S_1,\R))$ et ${h^{-1}}^* \circ \xi_\sigma^{-1}(H^1(S_2,\R))$ : en effet, soient $\omega_1 \in H^1(S_1,\R)$ et $\omega_2 \in H^1(S_2,\R)$. Alors
$$ \< *_X {h^{-1}}^* \circ \xi_\sigma^{-1}(\omega_1),{h^{-1}}^* \circ \xi_\sigma^{-1}(\omega_2)> = \Int({h^{-1}}^* \circ \xi_\sigma^{-1}(\omega_1),{h^{-1}}^* \circ \xi_\sigma^{-1}(\omega_2)) = 0,$$
donc $*_X{h^{-1}}^* \circ \xi_\sigma^{-1}(\omega_1)$ appartient à l'orthogonal de ${h^{-1}}^* \circ \xi_\sigma^{-1}(H^1(S_2,\R))$, qui est égal à ${h^{-1}}^* \circ \xi_\sigma^{-1}(H^1(S_1,\R))$.

\bigskip

La surface $X$ est conformément équivalente à une surface hyperbolique compacte connexe orientée $X'$ (marquée par $h' : S \ra X'$), privée de $q$ points (marqués par $P$). D'après l'invariance conforme (voir par exemple~\cite[Proposition~2.15, p.~27]{carron}), l'application qui à une $1$-forme différentielle harmonique sur $X'$ associe sa restriction à $X$ est un isomorphisme de $H^1(X',\R)$ dans $H^1_c(X,\R)$. Par ailleurs, l'étoile de Hodge ne dépend que de la structure complexe, donc l'isomorphisme entre $H^1(X',\R)$ et $H^1_c(X,\R)$ préserve les étoiles de Hodge. Ainsi nous savons que l'étoile de Hodge $*_{X'}$ préserve les sous-espaces ${h'^{-1}}^* \circ \xi_\sigma^{-1}(H^1(S_1,\R))$ et ${h'^{-1}}^* \circ \xi_\sigma^{-1}(H^1(S_2,\R))$ de $H^1(X',\R)$. Montrons que la fonction de longueur $\ell'=\psi'([X',h'])$ se scinde selon $\xi_\sigma$, avec $\psi' : \Tor(S,\emptyset) \ra {\R_+}^{\!H^1(S,\Z)}$.

\bigskip

Soient $\omega_1 \in H^1(S_1,\Z)$ et $\omega_2 \in H^1(S_2,\Z)$. Alors $\ell'(\xi_\sigma(\omega_1,\omega_2))^2$ est égal à
$$ \| {h'^{-1}}^* \circ \xi_\sigma(\omega_1,0)\|^2 + 2\<{h'^{-1}}^* \circ \xi_\sigma(\omega_1,0),{h'^{-1}}^* \circ \xi_\sigma(0,\omega_2)\> + \| {h'^{-1}}^* \circ \xi_\sigma(0,\omega_2)\|^2.$$
Or
$$ \<{h'^{-1}}^* \circ \xi_\sigma(\omega_1,0),{h'^{-1}}^* \circ \xi_\sigma(0,\omega_2)\> = \Int({h'^{-1}}^* \circ \xi_\sigma(\omega_1,0),*_{X'}{h'^{-1}}^* \circ \xi_\sigma(0,\omega_2))$$
est égal à $0$, car ${h^{-1}}^* \circ \xi_\sigma^{-1}(H^1(S_1,\R))$ et ${h^{-1}}^* \circ \xi_\sigma^{-1}(H^1(S_2,\R))$ sont stables par $*_{X'}$ et sont orthogonaux pour la forme d'intersection. Ainsi
$$ \ell'(\xi_\sigma(\omega_1,\omega_2))^2 = \ell'(\xi_\sigma(\omega_1,0))^2 + \ell'(\xi_\sigma(0,\omega_2))^2,$$
donc $\ell'$ se scinde selon $\xi_\sigma$.

\bigskip

Définissons alors l'élément $t$ de $\Sp_{2g}(\Z) = \Aut(H^1(S,\Z))$ par
\beq t : H^1(S,\Z) & \ra & H^1(S,\Z) \\
\omega & \mapsto & \xi_\sigma \circ \left( \left(\Id_{H^1(S_1,\Z)}\right) \oplus \left(-\Id_{H^1(S_2,\Z)}\right) \right) \circ {\xi_\sigma}^{-1}(\omega) \\
& & = \xi_\sigma((\xi_\sigma^{-1})_1(\omega),-(\xi_\sigma^{-1})_2(\omega)). \eeq 
Puisque la fonction de longueur $\ell'$ scinde selon $\xi_\sigma$, il est immédiat que pour tout $\omega \in H^1(S,\Z)$ nous avons $\ell'(t(\omega)) = \ell'(\omega)$.
L'action de $\Sp_{2g}(\Z)$ sur ${\R_+}^{\!H^1(S,\Z)}$ est par précomposition par l'adjoint, donc $t^* \cdot \ell' = \ell' \circ t = \ell'$. Ainsi $\psi(t^* \cdot [X',h']) = \psi([X',h'])$ dans ${\R_+}^{\!H^1(S,\Z)}$, donc $p'(t^* \cdot [X',h']) = p'([X',h'])$ dans $\mathcal{E}^b_{2g}$ par injectivité de l'application $\phi : \mathcal{E}^b_{2g} \ra \P({\R_+}^{\!H^1(S,\Z)})$. Par compacité de $X'$ on peut appliquer le théorème~\ref{thm:torelli}, donc on déduit que $t^* \cdot [X',h']$ est soit égal à $[X',h']$, soit à l'image de $[X',h']$ par l'involution hyperelliptique de l'espace de Torelli. Puisque l'involution hyperelliptique de $\Tor(S,\emptyset)$ est induite par l'élément $-\Id \in \Sp_{2g}(\Z)$, quitte à échanger les indices de $S_1$ et $S_2$, on peut supposer que $t^* \cdot [X',h'] = [X',h']$ dans $\Tor(S,\emptyset)$.

\bigskip

Considérons un difféomorphisme $\tau$ de $S$ préservant l'orientation tel que l'image de $\tau$ dans $\Sp_{2g}(\Z)$ soit égale à $t^*$, et tel que de plus $\tau \cdot [X',h'] = [X',h']$ dans l'espace de Teichmüller. \`{A} isotopie près, on peut supposer que le difféomorphisme $h \circ \tau \circ h^{-1}$ est une isométrie de $X$. Puisque $\tau^*$ préserve ${h'^{-1}}^* \circ \xi_\sigma^{-1}(H^1(S_1,\R))$ et ${h'^{-1}}^* \circ \xi_\sigma^{-1}(H^1(S_2,\R))$, on en déduit que $\tau$ préserve la classe d'homotopie de la courbe $\gamma$, or $h \circ \tau \circ h^{-1}$ est une isométrie donc stabilise la géodésique $h(\gamma)$.

\bigskip

Le difféomorphisme $\tau$ stabilise les deux composantes connexes $S_1 \bs P'_1$ et $S_2 \bs P'_2$ de $S \bs \gamma$, notons $\tau_1$ et $\tau_2$ les difféomorphismes de $S_1 \bs P'_1$ et $S_2 \bs P'_2$ induits. Puisque $P'_1$ et $P'_2$ sont des singletons, à isotopie près on peut supposer que $\tau_1$ et $\tau_2$ s'étendent en des difféomorphismes de $S_1$ et $S_2$ respectivement. Le difféomorphisme $\tau_1$ induit l'identité sur $H^1(S_1,\Z)$, donc sa classe d'isotopie appartient au groupe de Torelli $T(S_1,P_1)$. Or ce groupe est sans torsion, et comme $\tau$ est une involution nous avons $\tau_1^2=\tau_1$, donc $\tau_1$ appartient à $\Diff_0(S_1,P_1)$. Ainsi l'isométrie $h \circ \tau \circ h^{-1}$ restreinte à $X_1$ est isotope à l'identité, donc est égale à l'identité de $X_1$. Ainsi l'isométrie $h \circ \tau \circ h^{-1}$ est égale à l'identité de $X$. Or cette isométrie induit $-\Id$ sur $H^1_c(X_2,\Z)$ : c'est une contradiction.
\ep

\bibliographystyle{smfalpha_perso}
\bibliography{biblio}

\sign

\end{document}